\documentclass[11pt, letterpaper]{amsart}

\usepackage{amsthm,amsmath,amssymb,natbib}
\usepackage[normal]{subfigure}
\usepackage{graphicx, color}

\newtheorem{thm}{Theorem}[section]

\newtheorem{lem}[thm]{Lemma}

\theoremstyle{definition}

\theoremstyle{remark}
\newtheorem{rem}[thm]{Remark}
\newtheorem{exa}[thm]{Example}

\newcommand{\vs}{\vspace*{3mm}}

\def\ba{{\bf a}}
\def\bB{{\bf B}}

\def\bd{{\bf d}}

\def\bH{{\bf H}}

\def\bM{{\bf M}}

\def\bP{{\bf P}}

\def\bQ{{\bf Q}}

\def\bu{{\bf u}}

\def\bv{\mathbf v}

\def\bX{{\bf X}}
\def\bx{{\bf x}}
\def\bY{{\bf Y}}
\def\by{{\bf y}}
\def\bZ{{\bf Z}}

\def\boldeta{\mbox{\boldmath $\eta$}}
\def\bGamma{{\bf \Gamma}}
\def\bgamma{\mbox{\boldmath $\gamma$}}
\def\bLambda{{\bf \Lambda}}

\def\bmu{\mbox{\boldmath $\mu$}}

\def\bomega{\mbox{\boldmath $\omega$}}

\def\bSigma{{\bf \Sigma}}

\def\btau{\mbox{\boldmath $\tau$}}

\def\rd{{\rm d}}

\def\cK{{\mathcal K}}
\def\cA{{\mathcal A}}

\title[Saddlepoint methods for conditional expectations]{Saddlepoint methods for conditional expectations with applications to risk management}

\author[]{Sojung Kim}
\address[Sojung Kim]{Department of Mathematical Sciences, Korea Advanced Institute of Science and Technology, Daejeon 305-701, South Korea}
\email{sojungkim@kaist.ac.kr}
\author[]{Kyoung-Kuk Kim}
\address[Kyoung-Kuk Kim]{Department of Industrial and Systems Engineering, Korea Advanced Institute of Science and Technology, Daejeon 305-701, South Korea}
\email{kkim3128@kaist.ac.kr}

\date{July 15, 2015}

\begin{document}

\maketitle

\begin{abstract}
The paper derives saddlepoint expansions for conditional expectations in the form of $\mathsf{E}[\overline{X} | \overline{\bY} = \ba]$ and $\mathsf{E}[\overline{X} | \overline{\bY} \geq \ba]$ for the sample mean of a continuous random vector $(X, \bY^\top)$ whose joint moment generating function is available. Theses conditional expectations frequently appear in various applications, particularly in quantitative finance and risk management. Using the newly developed saddlepoint expansions, we propose fast and accurate methods to compute the sensitivities of risk measures such as value-at-risk and conditional value-at-risk, and the sensitivities of financial options with respect to a market parameter. Numerical studies are provided for the accuracy verification of the new approximations.

\vs
\end{abstract}

\baselineskip18pt

\section{Introduction}\label{Introduction}

 The saddlepoint method is one of the most important asymptotic approximations in statistics.
 It approximates a contour integral of Laplace type in the complex plane via the steepest descent method after a deformation of the original contour in such a way to contain the path of the steepest descent near the saddlepoint.  Since the development of saddlepoint approximations for the density of the sample mean of $n$ i.i.d. random variables by \cite{Daniels:54}, there have been numerous articles, treatises, and monographs on the topic. Their practical values have been particularly emphasized due to both high precision and simple explicit formulas.

 \cite{Barndorff-Cox:79} and \cite{Reid:88} initiated statistical applications of the saddlepoint method in inference such as approximating the densities of maximum likelihood estimators, likelihood ratio statistic or M-estimates. The widespread applicability in statistics also includes Bayesian analysis \citep{TIERNEY KADANE:86, Reid:03} and bootstrap inference \citep{BHW:92, Butler Bronson:02}.  Another important application is on financial option pricing and portfolio risk measurements in quantitative finance. From the opening paper of \cite{Rogers:99}, the saddlepoint method has been successfully applied in various contexts such as L\'evy processes \citep{CM:09}, affine jump-diffusion processes \citep{GK:09}, credit risk models \citep{Gordy:02} or value-at-risk \citep{Martin et al:01}, just to name a few. In such applications, one is usually concerned with obtaining approximate formulae for the density or tail probabilities of a target random variable.

 Relevant to this paper, pricing of collateralized debt obligations and computation of conditional value-at-risk requires evaluating the expectation in the form of $\mathsf{E}[Y {\bf 1}_{[Y  \geq a]}]$ for a random variable $Y$ and a constant $a$. Saddlepoint approximations to this expectation are derived in \cite{ Martin:06} and \cite{Huang Oost:11}. See Section~\ref{sec:pre2} for more details. Along the same line, the conditional expectations of the forms  $\mathsf{E}[X | Y = a]$ and $\mathsf{E}[X | Y \geq a]$ for a bivariate random vector $(X,Y)$ also appear in financial applications, but their saddlepoint approximations are not yet developed to the best of our knowledge.

 Let $(X, \bY^\top)$ be a continuous random vector where $X$ is a one-dimensional random variable and $\bY$ is a $d$-dimensional random vector. The objective of this paper is to derive saddlepoint expansions for conditional expectations in the form of $\mathsf{E}[\overline{X} | \overline{\bY} = \ba]$ and $\mathsf{E}[\overline{X} | \overline{\bY} \geq \ba]$ for the sample mean $\overline{X} = n^{-1} \sum_{i=1}^n X_i$ and $\overline{\bY} = n^{-1} \sum_{i=1}^n \bY_i$ with $\ba \in \mathbb{R}^d$. Here, the events $[\overline{\bY} = \ba]$ and $[\overline{\bY} \geq \ba]$ indicate the intersections of the respective univariate events. The derivation postulates the classical assumption of the existence of the joint density and the joint cumulant generating function $\mathcal{K}_{X, \bY}(\gamma, \boldeta)$ of $(X, \bY)$ which is analytic at the origin. We impose an additional assumption of an analytic property for the first derivative of the joint cumulant generating function with respect to the component of $X$ evaluated at zero,  $\mathcal{K}_{\gamma}(\eta) \triangleq {\partial}/{\partial \gamma} \{ \mathcal{K}_{X, \bY}(\gamma, \boldeta) \} |_{\gamma=0}$.

 Our first contribution is the derivation of saddlepoint approximations to the conditional expectations when $d = 1$ up to the order $O(n^{-2})$. As illustrated via several examples, the expansions are simple to apply and very accurate even for the case $n = 1$. The terms in the expansions only require the knowledge of the saddlepoint for the variable $Y$ and the derivatives of the cumulant generating function $\mathcal{K}_{Y}(\eta)$ of $Y$ and $\mathcal{K}_{\gamma}(\eta)$ evaluated at the saddlepoint.

 The second contribution is that the saddlepoint expansions for $d = 1$ are extended to the multivariate setting for $d \geq 2$. While the saddlepoint method for $\mathsf{E}[\overline{X} | \overline{\bY} = \ba]$ can be directly handled as in the case $d=1$, a major difficulty arises when deriving an expansion of $\mathsf{E}[\overline{X} | \overline{\bY} \geq \ba]$ due to the pole of the integrand. {To resolve this problem, we adopt the ideas presented in \cite{Kolassa:03} and \cite{Kolassa Li:10} where the authors study multivariate saddlepoint approximations.
 We decompose our target integrals into certain forms, for each of which the existing methods can be exploited.

 Last but not least, our saddlepoint approximations are demonstrated to be quite valuable in risk management. Either for portfolio risk measurements or hedging of financial contracts, it is important for a risk manager to know their sensitivities with respect to a specific parameter in order to make decisions in a responsive manner. Specifically in this work, we focus on the two widely popular risk measures, value-at-risk and conditional value-at-risk, and propose fast computational methods for their sensitivities by applying the newly developed saddlepoint expansions. Additionally, we show that sensitivities of an option based on multiple assets can be computed via the saddlepoint method. Numerical examples illustrate the effectiveness of our expansions in comparison with simulation based estimates.

 The rest of this paper is organized as follows. Section 2 first reviews classical saddlepoint approximations. Section 3 derives saddlepoint approximations to the target conditional expectations for $d=1$. The results in Section 3 are then extended to the multivariate setting in Section 4. Section 5 presents various applications in risk management with numerical studies. Finally, Section 6 concludes the paper.

\section{Preliminaries}\label{ch2}

\subsection{Classical saddlepoint approximation}\label{sec:Classical SPA}

 Let $\bY_1,\cdots, \bY_n$ be i.i.d. copies of a continuous random vector $\bY$ in $\mathbb{R}^d$ defined on a given probability space $(\Omega,\mathcal{F},\mathsf{P})$.
 We assume that $\bY$ has a bounded probability density function (PDF) and that its moment generating function (MGF) $m(\bgamma)$ exists for $\bgamma$ in some domain $\bGamma \subset \mathbb{R}^d$ containing an open neighborhood of the origin. The cumulant generating function (CGF) of $\bY$ is $\kappa(\bgamma) = \log m(\bgamma)$ defined in the same domain $\bGamma$.

 To describe classical saddlepoint techniques, we begin by recalling the inversion formula of the PDF and the tail probability of $\bY$: for $\by \in \mathbb{R}^d$,
    \begin{eqnarray}
        f_{\bY}(\by) &=& \left( \frac{1}{2\pi i} \right)^d \int_{\btau-i\infty}^{\btau+i\infty} \exp \left( \kappa(\bgamma) - \by^\top \bgamma \right) \rd \bgamma \mbox{ where } \btau \in \mathbb{R}^d; \label{multi-pdf-inversion} \\
        \mathsf{P}[\bY \geq \by] &=& \left(\frac{1}{2\pi i}\right)^d \int_{\btau-i\infty}^{\btau+i\infty} \frac{\exp \left( \kappa(\bgamma) - \by^\top \bgamma \right)}{\prod_{j=1}^d \gamma_j} \rd \bgamma \label{multi-cdf-inversion}
    \end{eqnarray}
 where $\gamma_j$ is the $j$-th component of $\bgamma$ and $\btau > {\bf 0 }\in \mathbb{R}^d$. We consider those values of $\by$ for which there exists the saddlepoint $\hat{\bgamma} = \hat{\bgamma}(\by)$ that solves the following saddlepoint equation
    \begin{equation*}\label{spaeqn}
        \kappa ' (\bgamma) = \by.
    \end{equation*}
 Throughout the paper, $f'(\bx)$ and $f''(\bx)$ of a multivariate function $f(\bx)$ denote its gradient and Hessian, respectively. The derivation of saddlepoint approximations first makes use of deformation of the original contour in the inversion formulas onto another contour containing the steepest descent curve that passes through the saddlepoint. After a suitable change of variable, asymptotic expansions of Laplace-type integrals are obtained with the help of Watson's lemma in \cite{Watson:48}.

 Let $\overline{\bY} = n^{-1} \sum_{i=1}^n \bY_i$ be the mean of $n$ i.i.d. observations.
 One classical saddlepoint approximation to the PDF of $\overline{Y}$ for $d=1$, known as Daniels' formula in  \cite{Daniels:54}, reads
     \begin{equation}\label{Daniels-SPA}
        f_{\overline{Y}}(y) = \sqrt{\frac{n}{2\pi \kappa''\left(\hat{\gamma}\right)}} e^{n\left[ \kappa\left(\hat{\gamma}\right) - y \hat{\gamma} \right]} \left[ 1 + \frac{1}{n} \left( \frac{\hat{\rho}_4}{8} - \frac{5 \hat{\rho}_3^2}{24}\right) + O\left(n^{-2}\right) \right]
    \end{equation}
 where $\hat{\rho}_r = \rho_r \left(\hat{\gamma}\right) = \kappa^{(r)}\left(\hat{\gamma}\right) / \kappa''\left(\hat{\gamma}\right)^{r/2}$ is the standardized cumulant of order $r$ evaluated at the saddlepoint $\hat{\gamma}$.

 For the tail probability of $\overline{Y}$, the Lugannani-Rice formula developed in \cite{Lugannani Rice:80} states
    \begin{eqnarray}\label{LR-SPA}
        \mathsf{P}[\overline{Y} \geq y] &=& \bar{\Phi}(\sqrt{n}\hat{\omega}) +  \frac{\phi(\sqrt{n}\hat{\omega})}{\sqrt{n}}\left[ \frac{1}{\hat{z}} - \frac{1}{\hat{\omega}} \right. \nonumber \\
        && \left. + \frac{1}{n} \left( \frac{1}{\hat{\omega}^3} - \frac{1}{\hat{z}^3} - \frac{\hat{\rho}_3}{2\hat{z}^2} + \left( \frac{\hat{\rho}_4}{8} - \frac{5 \hat{\rho}_3^2}{24}\right)\frac{1}{\hat{z}}\right) + O(n^{-2}) \right]
    \end{eqnarray}
 for $\hat\gamma$ away from zero where $\hat{\omega} = \textrm{sign}\left(\hat{\gamma}\right) \sqrt{2\left(y \hat{\gamma} - \kappa\left(\hat{\gamma}\right)\right)}$ and $\hat{z} = \hat{\gamma}\sqrt{\kappa''\left(\hat{\gamma}\right)}$. When $\hat{\gamma}$ is near zero, both $\hat{\omega}$ and $\hat{z}$ go to zero. Thus a different saddlepoint expansion should be employed in this case, for example, the formula (3.11) in \cite{Daniels:87}. The symbol $g(n) = O(n^\alpha)$ means that there exists a positive constant $C$ such that $|g(n)| \leq C n^\alpha$ as $n$ goes to infinity. The symbols $\phi(\cdot)$ and $\Phi(\cdot)$ denote the PDF and the cumulative distribution function (CDF) of a standard normal random variable, respectively. Lastly, $\bar{\Phi}(\cdot) = 1 - \Phi(\cdot)$.


 Such approximations for the PDF and the tail probability have their versions in the multivariate setting. A multivariate saddlepoint expansion of the PDF for a random vector $\overline{\bY}$ can be easily derived by extending Daniels' formula, and is presented as follows:
 \begin{equation}\label{multivariate-PDF-SPA}
        f_{\overline{\bY}}(\by) = \left(\frac{n}{2\pi}\right)^d \frac{\exp\left[ n \left( \kappa(\hat{\bgamma}) - \hat{\bgamma}^\top \by \right)\right]}{\sqrt{ \det \left[\kappa''\left(\hat{\bgamma}\right)\right]}} \left[ 1 + \frac{1}{n} \left( \frac{\hat{\varrho}_4 }{8} - \frac{\hat{\varrho}_{13}}{8} - \frac{\hat{\varrho}_{23}}{12} \right) + O(n^{-2}) \right]
    \end{equation}
  where the quantities $\hat{\varrho}_{4}$, $ \hat{\varrho}_{13}$, and $\hat{\varrho}_{23}$ are multivariate skewness and kurtosis, defined by
    \begin{eqnarray*}
        \hat{\varrho}_{4} &=& \sum_{i,j,p,l}\hat{\kappa}^{ijpl} \hat{\kappa}_{ij} \hat{\kappa}_{pl}, \\
        \hat{\varrho}_{13} &=& \sum_{i,j,p,l,m,o} \hat{\kappa}^{ijp} \hat{\kappa}^{lmo} \hat{\kappa}_{ij} \hat{\kappa}_{pl} \hat{\kappa}_{mo}, \mbox{ and }\\
        \hat{\varrho}_{23} &=& \sum_{i,j,p,l,m,o} \hat{\kappa}^{ijp} \hat{\kappa}^{lmo} \hat{\kappa}_{il} \hat{\kappa}_{jm} \hat{\kappa}_{po}.
    \end{eqnarray*}
 Here, the superscripted $\hat{\kappa}$ denotes the cumulants of the tilted distribution, that is, the derivatives of $\kappa(\bgamma) - \bgamma^\top \by$ evaluated at $\hat{\bgamma}$. For example, $\hat{\kappa}^{ijp} = \partial^3 \kappa(\bgamma) / \partial \gamma_i \partial \gamma_j \partial \gamma_p |_{\bgamma = \hat{\bgamma}} $. The subscripted  $\hat{\kappa}_{ij}$ refers to the $(i,j)$- entry of the inverse of the matrix formed by $\hat{\kappa}^{ij}$. The derivation of the terms is found in \cite{McCullagh:87}.

 On the other hand, the multivariate extension of saddlepoint expansions for the tail probability is somewhat difficult to achieve. Recently, \cite{Kolassa:03} and \cite{Kolassa Li:10} develop saddlepoint techniques to obtain an expansion up to the order $O(n^{-1})$; for a bivariate vector, see \cite{Wang:91}. Details are omitted here, but the key approaches of \cite{Kolassa:03} and \cite{Kolassa Li:10} appear in the multivariate version of our results in Section \ref{ch4}.

 For a detailed account of saddlepoint techniques, the reader is referred to \cite{Jensen:95}, \cite{Kolassa:06} or \cite{Butler:07}.

\subsection{Saddlepoint approximation to $\mathsf{E}[\overline{Y} | \overline{Y} \geq a]$}\label{sec:pre2}

 Interestingly, saddlepoint approximations to one special case of conditional expectation have been investigated, regarding the computation of conditional value-at-risk or also known as expected shortfall, a well-known risk measure defined as $\mathsf{E}[L | L \geq v_\alpha(L)]$ for a continuous random loss $L$ and value-at-risk $v_\alpha(L)$ of $L$ at level $\alpha$.

 When $L = \overline{Y}$ as in Section \ref{sec:Classical SPA}, one approach is to apply saddlepoint techniques to the integral
 $$
 \int_{-\infty}^a y f_{\overline{Y}}(y) \rd y.
 $$
 We first write $\mathsf{E}\left[\overline{Y} \ {\bf 1}_{[\overline{Y} \geq a]}\right]$ as $\mu - \int_{-\infty}^{a} y f_{\overline{Y}}(y) \rd y $, $\mu = \mathsf{E}[Y]$, and replace $f_{\overline{Y}}$ by Daniels' formula (\ref{Daniels-SPA}). And then an approximation to the integral of the form
    \begin{equation*}
        \sqrt{\frac{n}{2 \pi}} \int_{-\infty}^a e^{-n\zeta^2 / 2} \Psi_n(\zeta) \rd\zeta
    \end{equation*}
 for some function $\Psi_n$ can be employed from \cite{Temme:82}. This leads to the following formula which is also observed in \cite{Martin:06} up to the order $O\left(n^{-3/2}\right)$:
     \begin{eqnarray}\label{temme-SPA}
         \mathsf{E}\left[\overline{Y} \ {\bf 1}_{[\overline{Y} \geq a]}\right] &=& \mu \bar{\Phi}(\sqrt{n}\hat{\omega}) + \phi(\sqrt{n}\hat{\omega}) \frac{1}{\sqrt{n}} \left[\frac{a}{\hat{z}} -\frac{\mu}{\hat{\omega}} \right. \\
          && + \left. \frac{1}{n} \left( \frac{\mu}{\hat{\omega}^3} - \frac{a}{\hat{z}^3} - \frac{a\hat{\rho}_3}{2\hat{z}^2} + \frac{a}{\hat{z}}\left(\frac{\hat{\rho}_4}{8} - \frac{5 \hat{\rho}_3^2}{24}\right) + \frac{1}{\hat{\gamma}\hat{z}} \right) \right] + O\left(n^{-5/2}\right). \nonumber
     \end{eqnarray}

 Moreover, \cite{Butler Wood:04} obtain approximations to the MGF and its logarithmic derivatives of a truncated random variable $X_{(a,b)}$ with the density $f_X(x){\bf 1}_{(a,b)}(x) / (F_X(b) - F_X(a))$ for a distribution $F_X$ of $X$. Setting $b = \infty$ and $X = \overline{Y}$ and evaluating their approximation for the logarithmic derivative at zero produce another expansion:
     \begin{equation*}\label{butler-SPA}
         \mathsf{E}\left[\overline{Y} \ {\bf 1}_{[\overline{Y} \geq a]}\right] = \mu \bar{\Phi}(\sqrt{n}\hat{\omega}) + \phi(\sqrt{n}\hat{\omega}) \frac{1}{\sqrt{n}} \left[\frac{a}{\hat{z}} -\frac{\mu}{\hat{\omega}} + \frac{1}{n} \left( \frac{\mu - a}{\hat{\omega}^3} + \frac{1}{\hat{\gamma}\hat{z}} \right) \right] + O\left(n^{-5/2}\right).
     \end{equation*}
 \cite{BrodaPaolella:10} summarize the above mentioned methods in detail.
 
 \section{Saddlepoint approximation to conditional expectations}\label{ch3}

 Consider a continuous multi-dimensional random vector $(X, \bY^\top) \in \mathbb{R}^{d+1}$ where $X$ is a one-dimensional random variable and $\bY$ is a $d$-dimensional random vector.
 We define the multivariate MGF of $(X, \bY^\top)$ to be $\mathcal{M}_{X, \bY}(\gamma, \boldeta) = \mathsf{E}[\exp(\gamma X + \boldeta^\top \bY)]$ and the corresponding CGF to be $\mathcal{K}_{X, \mathbf{Y}}(\gamma, \boldeta) = \log \mathcal{M}_{X, \bY}(\gamma, \boldeta)$ for $\gamma \in \mathbb{R}$ and $\boldeta \in \mathbb{R}^d$. Classical assumptions are imposed: the joint PDF of $(X, \bY^\top)$ exists and the convergence domain of the CGF $\mathcal{K}_{X, \bY}(\gamma,\boldeta)$ contains an open neighborhood of the origin. The marginal CGFs of $X$ and $\bY$ are denoted by $\mathcal{K}_{X}(\gamma)$ and $\mathcal{K}_{\bY}(\boldeta)$, respectively.

 The goal of this section is to derive saddlepoint approximations to conditional expectations in the form of $\mathsf{E}[\overline{X} | \overline{\bY} = \ba]$ and $\mathsf{E}[\overline{X} | \overline{\bY} \geq \ba]$ for $\ba \in \mathbb{R}^d$ where $\overline{X} = n^{-1} \sum_{i=1}^n X_i$ and $\overline{\bY} = n^{-1} \sum_{i=1}^n \bY_i$ are the means of $n$ i.i.d. copies of $X$ and $\bY$, respectively. Thanks to the known formulas for PDFs and tail probabilities, the problem is reduced to utilizing saddlepoint techniques for $\mathsf{E}\left[\overline{X} \textbf{1}_{[\overline{\bY} = \ba]}\right]$ and $\mathsf{E}\left[\overline{X} \textbf{1}_{[\overline{\bY} \geq \ba]}\right]$.

 We first derive multivariate inversion formulas for $\mathsf{E}\left[X \textbf{1}_{[\bY = \ba]}\right]$ and $\mathsf{E}\left[X \textbf{1}_{[\bY \geq \ba]}\right]$ which resemble (\ref{multi-pdf-inversion}) and (\ref{multi-cdf-inversion}), respectively. We adopt the measure change approach of \cite{Huang Oost:11}.

 \begin{lem}\label{lem:inversion}
    For a continuous multivariate random vector $(X, \bY^\top) \in \mathbb{R}^{d+1}$, the following relations hold for $\btau$ in the domain of $\mathcal{K}_{\bY}$.
    \begin{equation}\label{lem:pdf-inversion}
        \mathsf{E}[X {\bf 1}_{[\bY = \ba]}] = \left( \frac{1}{2\pi i} \right)^d \int_{\btau -i\infty}^{\btau +i\infty} \left. \frac{\partial}{\partial \gamma} \mathcal{K}_{X, \bY}(\gamma, \boldeta) \right|_{\gamma=0} \exp \left( \mathcal{K}_{\bY}(\boldeta) - \ba^\top \boldeta \right) \rd\boldeta
    \end{equation}
    for $\btau \in \mathbb{R}^d$; and
    \begin{equation}\label{lem:cdf-inversion}
        \mathsf{E}[X {\bf 1}_{[\bY \geq \ba]}] = \left(\frac{1}{2\pi i}\right)^d \int_{\btau-i\infty}^{\btau + i\infty} \left. \frac{\partial}{\partial \gamma} \mathcal{K}_{X, \bY}(\gamma, \boldeta) \right|_{\gamma=0} \frac{\exp \left( \mathcal{K}_{\bY}(\boldeta) - \ba^\top \boldeta \right)}{\prod_{j=1}^d \eta_j} \rd \boldeta
    \end{equation}
    for $\btau > \mathbf{0}$
    where $\eta_j$ is the $j$-th component of $\boldeta$.
 \end{lem}
 \proof{ See Appendix \ref{app:proof of lemma1}.
 \endproof}

 In Sections \ref{sec:ch3-thm1} and \ref{sec:ch3-thm2}, we focus only on a bivariate random vector $(X,Y)$ where $d=1$ for its practical importance. In general, bivariate saddlepoint approximation requires to have a pair of saddlepoints that solve a system of saddlepoint equations, each of which depends on its respective variable. See, for example, \cite{DaYo:91}.
 However, in our development, only one saddlepoint of $\mathcal{K}_{Y}(\eta)$ is needed. Throughout the section, the saddlepoint $\hat{\eta} = \hat{\eta}(a)$ of $\mathcal{K}_{Y}(\eta)$ is assumed to exist as a solution of the saddlepoint equation
    \begin{equation}\label{spaeqn1}
        \frac{ \partial\mathcal{K}_{Y}}{\partial \eta}(\eta) = a.
    \end{equation}
 The conditions for the existence of a saddlepoint are discussed in Section 6 of \cite{Daniels:54}.

 \subsection{Saddlepoint approximation to $\mathsf{E}[\overline{X} | \overline{Y} = a]$}\label{sec:ch3-thm1}

 Before moving onto the derivation of an approximation to $\mathsf{E}\left[\overline{X} \textbf{1}_{[\overline{Y} = a]}\right]$, we shall present Watson's lemma which is the main technique to obtain an asymptotic expansion in powers of $n^{-1}$ in the classical approach. Our derivation relies on Watson's lemma applied to our new inversion formula in Lemma \ref{lem:inversion}. Here, its rescaled version is stated.

 \begin{lem}[Lemma 4.5.2 in \cite{Kolassa:06}]\label{Watsons}
    If $\vartheta(\omega)$ is analytic in a neighborhood of $\omega = \hat{\omega}$ containing the path $(-Ai + \hat{\omega}, Bi + \hat{\omega})$ with $0 < A, B \leq \infty$, then
    \begin{equation*}
        i^{-1} \left( \frac{n}{2\pi} \right)^{\frac{1}{2}} \int_{-Ai + \hat{\omega}}^{Bi + \hat{\omega}} \exp \left( - \frac{n}{2} (\omega - \hat{\omega})^2\right) \vartheta(\omega) \rd\omega = \sum_{j=0}^\infty \frac{(-1)^j \vartheta^{(2j)}(\hat{\omega})}{(2n)^j j!}
    \end{equation*}
    is an asymptotic expansion in powers of $n^{-1}$, provided the integral converges absolutely for some $n$.
 \end{lem}

 From the inversion formula (\ref{lem:pdf-inversion}) and the relations $\mathcal{K}_{\overline{X},\overline{Y}} (\gamma,\eta) = n \mathcal{K}_{X, Y}(\gamma/n, \eta/n)$ and $\mathcal{K}_{\overline{Y}}(\eta) = n \mathcal{K}_{Y}(\eta/n)$, the first target integral (\ref{lem:pdf-inversion}) is changed to
    \begin{equation}\label{inv1}
        \mathsf{E}\left[\overline{X} {\bf 1}_{[\overline{Y} = a]}\right] = \frac{n}{2\pi i} \int_{\tau-i\infty}^{\tau+i\infty} \left. \frac{\partial}{\partial \gamma} \mathcal{K}_{X, Y}(\gamma, \eta) \right|_{\gamma=0} \exp \left[ n \left( \mathcal{K}_{Y}(\eta) - a \eta \right)\right] \rd\eta
    \end{equation}
 for some $\tau\in \mathbb{R}$. For notational simplicity, we define
    \begin{equation*}
        \mathcal{K}_{\gamma}(\eta) \triangleq \left. \frac{\partial}{\partial \gamma} \mathcal{K}_{X, Y}(\gamma, \eta) \right|_{\gamma=0}.
    \end{equation*}
 We exploit the classical results to approximate (\ref{inv1}) but need to be careful when dealing with $\mathcal{K}_{\gamma}(\eta)$ in front of the exponential term.
 The next theorem is our first saddlepoint expansion for conditional expectation.

 \begin{thm}\label{thm1}
    Suppose that $\mathcal{K}_{\gamma}(\eta)$ is analytic in a neighborhood of $\hat{\eta}$.
    The conditional expectation $\mathsf{E}[\overline{X} | \overline{Y} = a] $ of a continuous bivariate random vector $(X,Y)$ can be approximated via saddlepoint techniques by
    \begin{eqnarray*}\label{thm1:spa}
        \mathsf{E}[\overline{X} | \overline{Y} = a] &=& \frac{1}{f_{\overline{Y}}(a)} \cdot \sqrt{\frac{n}{2\pi}} \frac{\exp \left[n\left(\mathcal{K}_{Y}(\hat{\eta}) - \hat{\eta}a \right) \right]}{\sqrt{\mathcal{K}''_{Y}(\hat{\eta})}} \left\{\mathcal{K}_{\gamma}(\hat{\eta}) + \frac{1}{n} \right. \nonumber\\
        && \times \left. \left[ \left( \frac{\hat{\rho}_4}{8} - \frac{5\hat{\rho}^2_3}{24} \right) \cdot \mathcal{K}_{\gamma}(\hat{\eta}) + \frac{\hat{\rho}_3}{2\sqrt{\mathcal{K}''_{Y}(\hat{\eta})}} \cdot \left.\frac{\partial}{\partial \eta }\mathcal{K}_{\gamma}(\eta)\right|_{\eta=\hat{\eta}} \right. \right. \\
        && \left.\left. - \frac{1}{2 \mathcal{K}''_{Y}(\hat{\eta})} \cdot \left.\frac{\partial^2}{\partial \eta^2}\mathcal{K}_{\gamma}(\eta)\right|_{\eta=\hat{\eta}} \right] + O(n^{-2}) \right\}, \nonumber
    \end{eqnarray*}
    where $\hat{\eta}$ is the saddlepoint that solves (\ref{spaeqn1}) and $\hat{\rho}_r = \mathcal{K}_{Y}^{(r)}\left(\hat{\eta}\right) / \mathcal{K}_{Y}''\left(\hat{\eta}\right)^{r/2}$ is the standardized cumulant of order $r$ evaluated at $\hat{\eta}$.

    Furthermore, if $f_{\overline{Y}}$ is also approximated by Daniel's formula (\ref{Daniels-SPA}), we have the following simple expansion:
    \begin{equation}\label{thm1:spa2}
        \mathsf{E}[\overline{X} | \overline{Y} = a] = \mathcal{K}_{\gamma}(\hat{\eta}) + \frac{\frac{\hat{\rho}_3}{2\sqrt{\mathcal{K}''_{Y}(\hat{\eta})}} \cdot \left.\frac{\partial}{\partial \eta }\mathcal{K}_{\gamma}(\eta)\right|_{\eta=\hat{\eta}} - \frac{1}{2 \mathcal{K}''_{Y}(\hat{\eta})} \cdot \left.\frac{\partial^2}{\partial \eta^2}\mathcal{K}_{\gamma}(\eta)\right|_{\eta=\hat{\eta}}}{ n+ \left( \frac{\hat{\rho}_4}{8} - \frac{5\hat{\rho}_3^2}{24} \right)} + O\left(n^{-2}\right).
    \end{equation}
 \end{thm}

 \proof{
 We integrate on the exactly same contour that is used in \cite{Daniels:54}. In Section 3 of \cite{Daniels:54}, the original path of integration is deformed into an equivalent path containing the steepest descent curve through the saddlepoint. On the steepest descent curve, the imaginary part of $\mathcal{K}_{Y}(\eta) - \eta a$ is a constant and its real part decreases fastest near $\hat{\eta}$. The contribution of the rest of the path to the target integral is negligible since some of them contribute a pure imaginary part and the others are bounded and converge to zero geometrically as $n$ goes to infinity.

 Rewrite (\ref{inv1}) using the closed curve theorem as
 \begin{eqnarray}\label{inv1-a}
    \mathsf{E}\left[\overline{X} {\bf 1}_{[\overline{Y} = a]}\right] &=& \frac{n}{2\pi i} \exp \left[n\left(\mathcal{K}_{Y}(\hat{\eta}) - \hat{\eta} a\right)\right] \times \\
    && \int_{\hat{\eta}-i\infty}^{\hat{\eta}+i\infty} \mathcal{K}_\gamma(\eta) \exp \Big[ n\Big(\mathcal{K}_{Y}(\eta) - \eta a - \mathcal{K}_{Y}(\hat{\eta}) + \hat{\eta} a\Big)\Big] \rd \eta. \nonumber
\end{eqnarray}
 The quantity in the exponent of the integrand, $\mathcal{K}_{Y}(\eta) - \eta a - \mathcal{K}_{Y}(\hat{\eta}) + \hat{\eta} a$, is an analytic function, and at $\hat{\eta}$ it is zero and has zero first derivative.

 Handling of the integrand in  \eqref{inv1-a} can be done via the classical approach well documented in, e.g., \cite{Kolassa:06}. Specifically, we make the same substitution (3.2) in \cite{Daniels:54} so that we have
 \begin{eqnarray*}
 \hat{\omega} &=& {\textrm{sign}(\hat{\eta})}\sqrt{2(\hat{\eta}a - \mathcal{K}_{Y}(\hat{\eta}) )},\\
 \omega(\eta) &=& \hat{\omega} + (\eta - \hat{\eta}) \sqrt{2 [\mathcal{K}_{Y}(\eta) - \eta a - \mathcal{K}_{Y}(\hat{\eta}) + \hat{\eta} a] / (\eta - \hat{\eta})^2}.
 \end{eqnarray*}
 Note that $\omega(\eta)$ is an analytic function of $\eta$ for $|\eta - \hat{\eta}| < \delta$ for some $\delta$, and by inverting the series of $\omega(\eta)$ we obtain an expansion of $\eta(\omega)$, the inverse of $\omega(\eta)$. Furthermore, it can be shown that
 \begin{equation}\label{dt_domega}
 \frac{\partial\eta}{\partial\omega} = \frac{1 - \frac{1}{3}\hat{\rho}_3(\omega(\eta) - \hat{\omega}) + (\frac{5}{24}\hat{\rho}_3^2 - \frac{1}{8}\hat{\rho}_4)(\omega(\eta) - \hat{\omega})^2 + O\left((\omega(\eta) - \hat{\omega})^3\right)}{\sqrt{\cK''_{Y}(\hat{\eta})}},
 \end{equation}
 whose verification is outlined in p.86 of \cite{Kolassa:06}.

 Then we re-parameterize (\ref{inv1-a}) in terms of $\omega$ as
 \begin{equation}\label{int_watson}
 \sqrt{\frac{n}{2\pi}} \exp \left[n\left(\mathcal{K}_{Y}(\hat{\eta}) - \hat{\eta} a\right)\right]  \times i^{-1}\sqrt{\frac{n}{2\pi}} \int_{\hat{\omega}-i\infty}^{\hat{\omega}+i\infty} \mathcal{K}_\gamma(\eta(\omega)) \exp \left[ \frac{n}{2} (\omega(\eta) - \hat{\omega})^2\right] \frac{\partial \eta}{\partial\omega} \rd\omega.
 \end{equation}
 Define
 $$
 \vartheta(\omega) \triangleq \mathcal{K}_\gamma(\eta(\omega)) \sqrt{\cK''_{Y}(\hat{\eta})} \frac{\partial \eta}{\partial\omega}.
 $$
 From the assumption on $\mathcal{K}_\gamma$ and the composition theorem of analytic functions, $\mathcal{K}_\gamma(\eta(\omega))$ has an expansion in a neighborhood of $\hat{\eta}$. And together with (\ref{dt_domega}), such an expansion leads us to conclude that $\vartheta(\omega)$ has a convergent series expansion in ascending powers of $\omega$. Then an asymptotic expansion of (\ref{int_watson}) is obtained directly from Lemma \ref{Watsons}, by inserting the expansion of $\vartheta(\omega)$ in (\ref{int_watson}) and integrating term-by-term:
 \begin{equation*}
 \mathsf{E}\left[\overline{X} {\bf 1}_{[\overline{Y} = a]}\right]=\sqrt{\frac{n}{2\pi}} \frac{\exp \left[n\left(\mathcal{K}_{Y}(\hat{\eta}) - \hat{\eta} a\right)\right]}{\sqrt{\cK''_{Y}(\hat{\eta})}} \left\{ \vartheta(\hat{\omega}) - \frac{1}{2n} \vartheta''(\hat{\omega}) + O(n^{-2})\right\}.
 \end{equation*}

 The first coefficient is $\vartheta(\hat{\omega}) = \mathcal{K}_{\gamma}(\hat{\eta})$. The second term is calculated from
 $$
 \vartheta''(\omega) = \sqrt{\cK''_{Y}(\hat{\eta})} \left\{ \frac{\partial^2}{\partial\eta^2}\mathcal{K}_\gamma(\eta) \left( \frac{\partial \eta}{\partial\omega}\right)^3 + 3 \frac{\partial}{\partial\eta}\mathcal{K}_\gamma(\eta) \frac{\partial \eta}{\partial\omega} \frac{\partial^2 \eta}{\partial\omega^2} + \mathcal{K}_\gamma(\eta(\omega)) \frac{\partial^3 \eta}{\partial\omega^3} \right\},
 $$
 differentiating (\ref{dt_domega}) with respect to $\omega$, and evaluating $\vartheta''(\omega)$ at $\hat{\omega}$. Detailed computations are omitted as they are straightforward.
 \endproof}

 In what follows, we illustrate some elementary examples in which the conditional expectation can be exactly calculated.

 \begin{exa}[Independent case]\label{thm1:ex1}
    When $X$ and $Y$ are independent, we have $\mathsf{E}[X | Y = a] = \mathsf{E}[X]$.
    Since $K_{X, Y}(\gamma,\eta) = \mathcal{K}_{X}(\gamma) + \mathcal{K}_{Y}(\eta)$, we have $\mathcal{K}_{\gamma}(\eta) = \mathcal{K}_{X}'(0)$ and $({\partial}/{\partial \eta})\mathcal{K}_{\gamma}(\eta) = ({\partial^2}/{\partial \eta^2})\mathcal{K}_{\gamma}(\eta) = 0.$
    Then (\ref{thm1:spa2}) turns out to be $\mathcal{K}'_{X}(0) = \mathsf{E}[X]$ which is exact.
 \end{exa}

 \begin{exa}[]\label{thm1:ex2}
    When $Y = X$, $\mathsf{E}[X | X = a] = a$. In that case, $K_{X, Y}(\gamma,\eta)  = \mathcal{K}_{X}(\gamma + \eta)$ and $\mathcal{K}_{\gamma}(\eta) = \mathcal{K}_{X}'(\eta)$. By computing $({\partial}/{\partial \eta})\mathcal{K}_{\gamma}(\eta) = \mathcal{K}_{X}''(\eta)$ and $({\partial^2}/{\partial \eta^2})\mathcal{K}_{\gamma}(\eta)= \mathcal{K}_{X}^{(3)}(\eta)$, the numerator of the second term in (\ref{thm1:spa2}) is zero; thus (\ref{thm1:spa2}) also results in $a$.
 \end{exa}

 \begin{exa}[Bivariate normal with correlation $\rho$]\label{thm1:ex3}
    Let $(X, Y)$ be a bivariate normal random variable, say $\mathcal{N}(\mu_1, \mu_2, \sigma_1^2, \sigma_2^2, \rho)$, with the CGF
    \begin{equation*}\label{cdf:bivariate normal}
        \cK_{X, Y}(\gamma,\eta) = \mu_1 \gamma + \mu_2 \eta + \frac{1}{2} \left( \sigma_1^2 \gamma^2 + 2 \rho \sigma_1 \sigma_2 \gamma \eta + \sigma_2^2 \eta^2 \right)
    \end{equation*}
    and correlation $\rho = {\sf Cov}(X,Y)/\sigma_1 \sigma_2$. Note that $(\overline{X},\overline{Y}) \sim \mathcal{N}(\mu_1, \mu_2, \sigma_1^2/n, \sigma_2^2/n, \rho)$.
    Thus,
    \begin{equation*}
        \mathsf{E}[\overline{X} | \overline{Y}=a] = \mu_1+ \rho \frac{\sigma_1}{\sigma_2} (a - \mu_2).
    \end{equation*}
    On the other hand, $\mathcal{K}_{\gamma}(\eta) = \mu_1 + \rho \sigma_1 \sigma_2 \eta$, $(\partial/\partial \eta)\mathcal{K}_{\gamma}(\eta) = \rho \sigma_1 \sigma_2$ and $(\partial^2/\partial \eta^2)\mathcal{K}_{\gamma}(\eta) = 0$.
    The saddlepoint $\hat{\eta}(a)$ is $\hat{\eta} = (a - \mu_2)/\sigma_2^2$ and the 3rd order standardized cumulant $\hat{\rho}_3$ is zero.
    Therefore, (\ref{thm1:spa2}) yields the exact result as $\mathsf{E}[\overline{X} | \overline{Y}=a] = \mathcal{K}_{\gamma}(\hat{\eta})$.
 \end{exa}

 \subsection{Saddlepoint approximation to $\mathsf{E}[\overline{X} | \overline{Y} \geq a]$}\label{sec:ch3-thm2}

 Under the setting of Section \ref{sec:ch3-thm1}, the second target integral can be rewritten by the inversion formula (\ref{lem:cdf-inversion}) as
    \begin{equation}\label{inv2}
        \mathsf{E}\left[\overline{X} {\bf 1}_{[\overline{Y} \geq a]}\right] = \left(\frac{1}{2\pi i}\right) \int_{\tau-i\infty}^{\tau+i\infty} \mathcal{K}_\gamma(\eta) \frac{\exp \left[ n \left( \mathcal{K}_{Y}(\eta) - a \eta \right)\right]}{\eta} \rd \eta
    \end{equation}
    for $\tau > 0$.
 Following the approach in \cite{Martin:06}, we divide the singularity in the integrand as
    \begin{equation*}
        \frac{\mathcal{K}_\gamma(\eta)}{\eta} = \frac{\mathcal{K}_\gamma(0)}{\eta} + \frac{\mathcal{K}_\gamma(\eta) - \mathcal{K}_\gamma(0)}{\eta}.
    \end{equation*}
 Then, (\ref{inv2}) becomes the sum of two tractable parts, namely, for $\tau > 0$
    \begin{equation}\label{inv2-2}
        \mathsf{E}\left[\overline{X} {\bf 1}_{[\overline{Y} \geq a]}\right] = \mathcal{K}_\gamma(0) \cdot \mathsf{P}[\overline{Y} \geq a] +
        \frac{1}{2\pi i} \int_{\tau-i\infty}^{\tau+i\infty} \frac{\mathcal{K}_\gamma(\eta) - \mathcal{K}_\gamma(0)}{\ \eta} \exp \left[ n \left( \mathcal{K}_{Y}(\eta) - a \eta \right)\right] d \eta.
    \end{equation}
 The second complex integral is treated in the similar fashion as in Theorem \ref{thm1}, using Lemma \ref{Watsons}.

 \begin{thm}\label{thm2}
    Suppose that $\mathcal{K}_{\gamma}(\eta)$ is analytic in a neighborhood of $\hat{\eta}$ and that $Y$ is continuous at $a$.
    The conditional expectation $\mathsf{E}[\overline{X} | \overline{Y} \geq a] $ of a continuous bivariate random vector $(X,Y)$ can be approximated via saddlepoint techniques by
    \begin{eqnarray}\label{thm2:spa}
        \lefteqn{\mathsf{E}[\overline{X} | \overline{Y} \geq a] = \mathsf{E}[X] \ + \frac{1}{\mathsf{P}[\overline{Y} \geq a]} \frac{1}{\sqrt{2\pi n}} \exp \left[n\left(\mathcal{K}_{Y}(\hat{\eta}) - \hat{\eta}a \right) \right]}\\
         &\times&  \left\{ \frac{\mathcal{K}_\gamma(\hat{\eta}) - \mathcal{K}_\gamma(0)}{\hat{z}} \right.
         +  \frac{1}{n} \left[ \frac{\mathcal{K}_\gamma(\hat{\eta}) - \mathcal{K}_\gamma(0)}{\hat{z}}
        \left( \frac{\hat{\rho}_4}{8} - \frac{5\hat{\rho}^2_3}{24} - \frac{\hat{\rho}_3}{2\hat{z}} - \frac{1}{\hat{z}^2}\right) \right. \nonumber \\
        && + \frac{1}{\hat{z} \sqrt{\mathcal{K}_{Y}''(\hat{\eta})}} \left( \frac{\hat{\rho}_3}{2} + \frac{1}{\hat{z}} \right) \cdot \left.\frac{\partial}{\partial \eta} \mathcal{K}_{\gamma}(\eta)\right|_{\eta=\hat{\eta}}  - \left.\left. \frac{1}{2 \hat{z} {\mathcal{K}_{Y}''(\hat{\eta})}} \cdot \left.\frac{\partial^2}{\partial \eta^2}\mathcal{K}_{\gamma}(\eta)\right|_{\eta=\hat{\eta}} \right] + O(n^{-2}) \right\}, \nonumber
    \end{eqnarray}
    where $\hat{\eta}$ solves (\ref{spaeqn1}), $\hat{z} = \hat{\eta}\sqrt{\mathcal{K}_{Y}''(\hat{\eta})}$, and $\hat{\rho}_r = \mathcal{K}_{Y}^{(r)}\left(\hat{\eta}\right) / \mathcal{K}_{Y}''\left(\hat{\eta}\right)^{r/2}$.

    When $\hat{\eta} = 0$, we have an expansion
    \begin{equation*}\label{thm2:spa2}
        \mathsf{E}[\overline{X} | \overline{Y} \geq a ] = \mathsf{E}[X] + \frac{1}{\sqrt{2\pi n \mathcal{K}''_{Y}(0)} \cdot \mathsf{P}[\overline{Y} \geq a]} \cdot   \left.\frac{\partial}{\partial \eta} \mathcal{K}_{\gamma}(\eta)\right|_{\eta=0} + O\left(n^{-3/2}\right).
    \end{equation*}
 \end{thm}

 \proof{
 Let $a = \mathcal{K}_{Y}'(\hat{\eta})$ and first we suppose that $\hat{\eta}> 0$, or equivalently $\mathsf{E}[Y] > a$.  Again, we only focus on the integration on the steepest descent curve and take the new variables $\omega$ and $\hat{\omega}$ as in the proof of Theorem \ref{thm1}.
 To expand $(1/\eta)(\partial \eta / \partial \omega)$, we closely follow the approach in p.92 of \cite{Kolassa:06}. First, we integrate the expansion in (\ref{dt_domega}) to obtain
 \begin{equation}\label{d_omega}
    \eta =  \hat{\eta} + \frac{1}{\sqrt{\mathcal{K}_{Y}''(\hat{\eta})}} \left[ (\omega - \hat{\omega}) - \frac{1}{6}\hat{\rho}_3 (\omega - \hat{\omega})^2 + \left( \frac{5}{72}\hat{\rho}_3^2 - \frac{1}{24}\hat{\rho}_4\right) (\omega - \hat{\omega})^3 + \ O((\omega - \hat{\omega})^4)\right].
 \end{equation}
 Then, dividing (\ref{dt_domega}) by (\ref{d_omega}) yields
 \begin{eqnarray}\label{tdt_domega}
    \frac{1}{\eta} \frac{\partial \eta}{\partial \omega} &=& \left[ 1 - \left(\frac{\hat{\rho}_3}{3} + \frac{1}{\hat{z}} \right)(\omega - \hat{\omega}) + \left(\frac{5}{24}\hat{\rho}_3^2 - \frac{1}{8}\hat{\rho}_4 + \frac{\hat{\rho}_3}{2}\frac{1}{\hat{z}} + \frac{1}{\hat{z}^2} \right)(\omega - \hat{\omega})^2 \right. \nonumber \\
    && + \ O\left((\omega - \hat{\omega})^3\right) \Big] \Big/ \hat{z}
 \end{eqnarray}
 where $\hat{z} =\hat{\eta} \sqrt{\mathcal{K}_{Y}(\hat{\eta})}$. Note that the coefficients of the odd order terms of $\omega - \hat{\omega}$ should be determined since it does not disappear in our derivation, whereas they are removed in the classical approach. See (101) of \cite{Kolassa:06}.

 Define
 \begin{equation*}\label{func_vt}
    \vartheta(\omega) \triangleq \left(\mathcal{K}_\gamma(\eta) - \mathcal{K}_\gamma(0)\right) \frac{1}{\eta}\frac{\partial \eta}{\partial\omega}
 \end{equation*}
 whose convergent series exists at $\hat{\omega}$ by (\ref{tdt_domega}) and the analytic property of  $\mathcal{K}_\gamma(\eta)$.
 Then the second term in (\ref{inv2-2}) becomes
 \begin{eqnarray}\label{int_es2}
    && \sqrt{\frac{1}{2\pi n}} \exp \left[n(\mathcal{K}_{Y}(\hat{\eta}) - \hat{\eta} a) \right] \cdot i^{-1} \sqrt{\frac{n}{2\pi}}\int_{\hat{\omega}-i\infty}^{\hat{\omega}+i\infty} \exp \left[\frac{n}{2}(\omega-\hat{\omega})^2 \right] \vartheta(\omega)\rd\omega \nonumber \\
    &=& \sqrt{\frac{1}{2\pi n}}  \exp \left[n(\mathcal{K}_{Y}(\hat{\eta}) - \hat{\eta} a) \right] \sum_{j=0}^{\infty} \frac{(-1)^j \vartheta^{(2j)}(\hat{\omega})}{(2n)^j j!}
 \end{eqnarray}
 by Watson's lemma.
 The coefficients in the expansion (\ref{int_es2}) are calculated by expanding $\vartheta(\omega)$ about $\hat{\omega}$.
 By combining (\ref{dt_domega}), (\ref{d_omega}) and (\ref{tdt_domega}), and taking their derivatives, we compute $\vartheta(\hat{\omega}) = (\mathcal{K}_\gamma(\hat{\eta}) - \mathcal{K}_\gamma(0))/{\hat{z}}$, and
 \begin{eqnarray*}
    \vartheta''(\hat{\omega}) &=&  2 \frac{\mathcal{K}_\gamma(\hat{\eta}) - \mathcal{K}_\gamma(0)}{\hat{z}}\left(\frac{5}{24}\hat{\rho}_3^2 - \frac{1}{8}\hat{\rho}_4 + \frac{\hat{\rho}_3}{2\hat{z}} + \frac{1}{\hat{z}^2} \right) \\
    && - \frac{1}{\hat{z} \sqrt{\mathcal{K}_{Y}''(\hat{\eta})}} \left( \hat{\rho}_3 + \frac{2}{\hat{z}} \right) \cdot \left. \frac{\partial}{\partial \eta} \mathcal{K}_\gamma(\eta) \right|_{\eta = \hat{\eta}} + \frac{1}{\hat{z} \mathcal{K}_{Y}''(\hat{\eta})} \cdot \left.\frac{\partial^2}{\partial \eta^2} \mathcal{K}_\gamma(\eta) \right|_{\eta = \hat{\eta}}.
 \end{eqnarray*}
 The desired result is then immediate.

 Now suppose that $\hat\eta < 0$. We set $Z = - Y$ and observe that
    \begin{equation*}
        \mathsf{E}\left[\overline{X} {\bf 1}_{[\overline{Y} \geq a]}\right] = \mathsf{E}\left[\overline{X} {\bf 1}_{[\overline{Z} \leq -a]}\right] = \mathsf{E}[X] - \mathsf{E}\left[\overline{X} {\bf 1}_{[\overline{Z} \geq -a]}\right].
    \end{equation*}
 For the second term on the right hand side, the saddlepoint that satisfies $\cK_Z'(\cdot) = -a$ is $-\hat\eta > 0$. Working with the CGF of $(X,Z)$ and transforming back to $Y$, an expansion for $\hat{\eta} < 0$ can be found. And the final formula turns out to be the same formula as (\ref{thm2:spa}).

 When $\hat{\eta}=0$, equivalently  $\hat{\omega}=0$,  $\lim_{\omega \rightarrow 0} \eta(\omega) = 0$ and
 $$ \lim_{\omega \rightarrow 0} \frac{\mathcal{K}_\gamma(\eta) - \mathcal{K}_\gamma(0)}{\eta} \frac{\partial\eta}{\partial\omega} = \frac{1}{ \sqrt{ \mathcal{K}_{Y}''(0)}} \cdot \left. \frac{\partial}{\partial \eta} \mathcal{K}_\gamma(\eta) \right|_{\eta = 0}.$$
 Thus, $\vartheta(\omega)$ is analytic at $\omega = 0$. This yields the following approximation to (\ref{int_es2}) for $\hat{\eta}=0$ by applying Watson's lemma centered at $\hat{\omega}=0$:
 \begin{equation*}
 \frac{1}{\sqrt{2\pi n} \cdot \sqrt{ \mathcal{K}_{Y}''(0)}} \cdot  \left. \frac{\partial}{\partial \eta} \mathcal{K}_\gamma(\eta) \right|_{\eta = 0} + O\left(n^{-3/2}\right).
 \end{equation*}
 \endproof}

 \begin{rem}
    The saddlepoint approximation to the lower-tail expectation $\mathsf{E}\left[\overline{X}{\bf 1}_{[\overline{Y} \leq a]}\right]$ can be obtained simply by considering $\mathsf{E}[X] - \mathsf{E}\left[\overline{X} {\bf 1}_{[\overline{Y} \geq a]}\right]$ and by using $\eqref{thm2:spa}$ for the second term. Alternatively, we can obtain an approximation to the integral directly by applying (\ref{thm2:spa}) by replacing $Y$ with $-Y$. In either case, the resulting formula is the same.
 \end{rem}

 \begin{exa}[Bivariate normal with correlation $\rho$]\label{thm2:ex1}
    Consider Example \ref{thm1:ex3} where $(X, Y) \sim \mathcal{N}(\mu_1, \mu_2, \sigma_1^2, \sigma_2^2, \rho)$.
    Evaluating \eqref{inv2-2} gives us
    \begin{eqnarray*}
    {\sf E}\left[\overline{X} |\overline{Y} \geq a \right] &=& \mu_1  + \rho\sigma_1\sigma_2 \frac{1}{2\pi i}\int_{\tau - i \infty}^{\tau + i\infty}\exp\left[n(\cK_Y(\eta) - a\eta)\right] \rd \eta / {\sf P}[\overline Y \geq a]\\
    &=& \mu_1  + \frac{\rho\sigma_1\sigma_2}{n} \phi_{\overline{Y}}(a) / {\sf P}[\overline Y \geq a]
    \end{eqnarray*}
    where $\phi_{\overline{Y}}$ is the PDF of $\overline{Y}$. On the other hand, it is easy to check that (\ref{thm2:spa}) yields the same value by
    \begin{equation*}
        \mathsf{E}[\overline{X} | \overline{Y} \geq a ] = \mu_1 + \frac{\rho \sigma_1}{\sqrt{n}} \phi(\sqrt{n}\hat{\omega}) / \mathsf{P}[\overline{Y} \geq a]
    \end{equation*}
    where
    $\hat{\omega} = {\textrm{sign}(\hat{\eta})}\sqrt{2(\hat{\eta}a - \mathcal{K}_{Y}(\hat{\eta}) )} = \textrm{sign}(a-\mu_2)(a-\mu_2)/\sigma_2$.
 \end{exa}

 \begin{rem}\label{rem:integralspa}
 By approximating $\mathsf{P}[\overline{Y} \geq a]$ with the Lugannani-Rice formula (\ref{LR-SPA}), the expansion \eqref{thm2:spa} for $\mathsf{E}\left[\overline{X} {\bf 1}_{[\overline{Y} \geq a]}\right] $ is reduced to
 \begin{eqnarray}\label{rem:spa}
      \lefteqn{\mu \bar{\Phi}(\sqrt{n}\hat{\omega}) + \phi(\sqrt{n}\hat{\omega}) \frac{1}{\sqrt{n}} \left[ \frac{\mathcal{K}_\gamma(\hat{\eta})}{\hat{z}} - \frac{\mu}{\hat{\omega}} + \frac{1}{n} \left( \frac{\mathcal{K}_\gamma(\hat{\eta})}{\hat{z}}
    \left( \frac{\hat{\rho}_4}{8} - \frac{5\hat{\rho}^2_3}{24} - \frac{\hat{\rho}_3}{2\hat{z}} - \frac{1}{\hat{z}^2}\right) + \frac{\mu}{\hat{\omega}^3} \right.\right.}  \nonumber \\
    && \left. \left. +  \frac{1}{\hat{z} \sqrt{\mathcal{K}_{Y}''(\hat{\eta})}} \left( \frac{\hat{\rho}_3}{2} + \frac{1}{\hat{z}} \right) \cdot \left.\frac{\partial}{\partial \eta} \mathcal{K}_{\gamma}(\eta)\right|_{\eta=\hat{\eta}}  - \frac{1}{2 \hat{z} {\mathcal{K}_{Y}''(\hat{\eta})}} \cdot \left.\frac{\partial^2}{\partial \eta^2}\mathcal{K}_{\gamma}(\eta)\right|_{\eta=\hat{\eta}} \right) \right]
 \end{eqnarray}
 where $\mu = \mathsf{E}[X]$.
 When $X = Y$, it becomes exactly the same as (\ref{temme-SPA}).
 \end{rem}

 Discussions about the accuracy of the expansions in Theorems \ref{thm1} and \ref{thm2} are deferred to Section \ref{ch5} where numerical studies are provided in the context of risk management.
 
 \section{Multivariate extension}\label{ch4}

 In this section, we consider the case $d \geq 2$. The saddlepoint $\hat{\boldeta} = \hat{\boldeta}(\ba)$ of $\mathcal{K}_{\bY}(\boldeta)$ is assumed to exist as the solution to the system of saddlepoint equations
    \begin{equation}\label{spaeqn-multi}
        \frac{ \partial\mathcal{K}_{\bY}}{\partial \eta_i}(\boldeta) = a_i,
    \end{equation}
 where $\ba = (a_1, \cdots, a_d)$ and $\boldeta = (\eta_1, \cdots, \eta_d)$ for $i = 1, \cdots, d$.
 As before, define
    \begin{equation*}
        \mathcal{K}_{\gamma}(\boldeta) \triangleq \left. \frac{\partial}{\partial \gamma} \mathcal{K}_{X, \bY}(\gamma, \boldeta) \right|_{\gamma=0}.
    \end{equation*}

\subsection{Extension of Theorem~\ref{thm1}}

 Finding an analog of Theorem \ref{thm1} for the case $d \geq 2$ raises no additional difficulty because we can utilize a multivariate version of Watson's lemma, which is also useful when deriving multivariate saddlepoint approximations to multivariate PDFs.
 To be specific, we take a differentiable function $\bomega(\boldeta)$ via the change of variable
    \begin{equation}\label{def:w1}
        \frac{1}{2}(\bomega - \hat{\bomega})^\top(\bomega- \hat{\bomega}) = \mathcal{K}_{\mathbf{Y}}(\boldeta) - \boldeta^\top \mathbf{a} -  \mathcal{K}_{\mathbf{Y}}(\hat{\boldeta}) + \hat{\boldeta}^\top \mathbf{a}
    \end{equation}
 which is employed in \cite{Kolassa:96}.
 This function is proved to be analytic for $\bomega$ in a neighborhood of $\hat{\bomega}$, and the detailed construction will be given for $d=2$ in the next subsection.
 Using the change of variable \eqref{def:w1} in (\ref{lem:cdf-inversion}) with $(\overline{X}, \overline{\bY})$, and applying multivariate Watson's lemma \ref{multi-Watsons} with a particular care for $ \mathcal{K}_{\gamma}(\boldeta)$, we arrive at the following result.

 \begin{thm}\label{cor:e-thm1}
    Let $\hat{\boldeta}$ be a solution to the saddlepoint equation (\ref{spaeqn-multi}) and suppose that $\mathcal{K}_\gamma(\boldeta)$ is analytic in a neighborhood of $\hat{\boldeta}$.
    The conditional expectation $\mathsf{E}[\overline{X} | \overline{\bY} = \ba] $ of a continuous random vector $(X,\bY)$ can be approximated via saddlepoint techniques by
    \begin{eqnarray*}\label{cor:spa}
        \mathsf{E}[\overline{X} | \overline{\bY} = \ba] &=& \frac{1}{f_{\overline{\bY}}(\ba)} \cdot \left(\frac{n}{2\pi}\right)^{d/2} \frac{\exp \left[n\left(\mathcal{K}_{\bY}(\hat{\boldeta}) - \hat{\boldeta}^\top \ba\right)\right]}{\sqrt{\det\left[\mathcal{K}''_{\bY}(\hat{\boldeta})\right]}} \nonumber \\
        &&\times  \left\{\mathcal{K}_\gamma(\hat{\boldeta}) + \frac{1}{2n} \left[ \mathcal{K}_\gamma(\hat{\boldeta}) \cdot \beta(\hat{\boldeta}) + \sum_{i} \left. \frac{\partial}{\partial \eta_i} \mathcal{K}_\gamma(\boldeta)\right|_{\boldeta = \hat{\boldeta}} \beta_{i}(\hat{\boldeta}) \right.\right. \\
        &&  \left.\left. + \sum_{i,j} \left.\frac{\partial^2}{\partial \eta_i \partial \eta_j} \mathcal{K}_\gamma(\boldeta) \right|_{\boldeta = \hat{\boldeta}} \beta_{i,j}(\hat{\boldeta})\right] + O\left(n^{-2}\right)\right\}. \nonumber
    \end{eqnarray*}
    The coefficients $\beta, \beta_i$, and $\beta_{i,j}$ evaluated at $\hat{\boldeta}$ satisfy
    \begin{align*}
    & \beta(\hat{\boldeta}) = -  \sum_{k=1}^d \frac{\partial^2}{\partial \omega_k^2}\left| \frac{\partial \boldeta}{\partial \bomega} \right|_{\bomega = \hat{\bomega}} \sqrt{\det\left[\mathcal{K}''_{\bY}(\hat{\boldeta})\right]};\\
    & \beta_{i}(\hat{\boldeta}) = - \sum_{k=1}^d \left\{ \frac{\partial^2 \eta_i}{\partial \omega_k^2} (\hat{\bomega}) + 2 \frac{\partial \eta_i}{\partial \omega_k}(\hat{\bomega}) \cdot \frac{\partial}{\partial \omega_k}\left| \frac{\partial \boldeta}{\partial \bomega} \right|_{\bomega = \hat{\bomega}} \cdot \sqrt{\det\left[\mathcal{K}''_{\bY}(\hat{\boldeta})\right]} \right\}; \mbox{ and }\\
    & \beta_{i,j}(\hat{\boldeta}) = - \sum_{k=1}^d \frac{\partial \eta_i}{\partial \omega_k}(\hat{\bomega})\frac{\partial \eta_j}{\partial \omega_k}(\hat{\bomega})  ,
    \end{align*}
    respectively.

    Furthermore, if $f_{\overline{\bY}}(\ba)$ is also approximated by Daniel's formula (\ref{multivariate-PDF-SPA}), we have the following simple expansion:
    \begin{eqnarray*}\label{cor:spa2}
        \lefteqn{ \mathsf{E}[\overline{X} | \overline{\bY} = \ba]      = \mathcal{K}_{\gamma}(\hat{\boldeta}) } && \\
         && + \frac{\sum_{i} \left. \frac{\partial}{\partial \eta_i} \mathcal{K}_\gamma(\boldeta)\right|_{\boldeta = \hat{\boldeta}} \beta_{i}(\hat{\boldeta}) + \sum_{i,j} \left.\frac{\partial^2}{\partial \eta_i \partial \eta_j} \mathcal{K}_\gamma(\boldeta) \right|_{\boldeta = \hat{\boldeta}} \beta_{i,j}(\hat{\boldeta})}{ 2n + \beta(\hat{\boldeta})} + O\left(n^{-2}\right).
    \end{eqnarray*}
 \end{thm}

 \proof{ See Appendix \ref{app:proof of cor}.
 \endproof}

 On the other hand, a major concern arises when deriving the extension of Theorem \ref{thm2}.
 Due to the factor in the denominator of (\ref{multi-cdf-inversion}) which is apparently not a simple pole, multivariate saddlepoint approximation to the tail probability is difficult to compute. Among various methods to tackle the problem, \cite{Kolassa Li:10} suggest an approach to extend the method of \cite{Lugannani Rice:80} to the multivariate case. The authors obtain a tractable formula up to the relative order $O\left(n^{-1}\right)$. We essentially adopt their framework but particular attention should be paid to the multiplying factor $\mathcal{K}_\gamma(\boldeta)$ in computing $\mathsf{E}\left[\overline{X} {\bf 1}_{[\overline{\bY} \geq \ba]}\right]$. Under a suitable assumption on $\cK_\gamma(\boldeta)$, we decompose $\cK_\gamma(\boldeta)$ in such a way that each corresponding integral can be approximated separately. In the next subsection, the extension of Theorem \ref{thm2} is stated for the case $d=2$ for an illustration and practical usefulness. The entire idea is still applicable when $d > 2$ but it is computationally heavy.

 \subsection{Extension of Theorem~\ref{thm2}}\label{sec:ch4-thm2}

 With $\bY \in \mathbb{R}^2$, the inversion formula is written as
    \begin{equation}\label{e-inv2}
        \mathsf{E}\left[\overline{X} {\bf 1}_{[\overline{\bY} \geq \ba]}\right] = \left(\frac{1}{2\pi i}\right)^2 \int_{\btau-i\infty}^{\btau+i\infty} \mathcal{K}_\gamma(\eta_1, \eta_2) \frac{\exp \left[ n \left( \mathcal{K}_{\bY}(\eta_1, \eta_2) - \eta_1 a_1 - \eta_2a_2 \right)\right]}{\eta_1 \eta_2} \rd \boldeta,
    \end{equation}
    for $\btau > \mathbf{0}$.
 In order to identify the pole in the integrand of (\ref{e-inv2}), we adopt the following explicit functions constructed in \cite{Kolassa Li:10}.

 Define $\tilde{\eta}_2(\eta_1)$ as the minimizer of $\mathcal{K}_{\bY}(\eta_1, \eta_2) - \eta_1 a_1 - \eta_2 a_2$ when the first component $\eta_1$ is fixed, i.e.,
 \begin{equation*}
    \tilde{\eta}_2(\eta_1) = \arg\min_{\eta_2} \left\{ \mathcal{K}_{\bY}(\eta_1, \eta_2) - \eta_1 a_1 - \eta_2 a_2\right\}.
 \end{equation*}
 The analytic function $\bomega(\boldeta)$ satisfying (\ref{def:w1}) is further specified as
    \begin{eqnarray*}\label{def:w3}
        -\frac{1}{2} \hat{\omega}_1^2 &=& \mathcal{K}_{\bY}(\hat{\eta}_1, \hat{\eta}_2) - \hat{\eta}_1 a_1 - \hat{\eta}_2 a_2 - \left( \mathcal{K}_{\bY}(0, \tilde{\eta}_2(0)) - \tilde{\eta}_2(0) a_2 \right),\\
        -\frac{1}{2} (\omega_1 - \hat{\omega}_1)^2 &=& \mathcal{K}_{\bY}(\hat{\eta}_1, \hat{\eta}_2) - \hat{\eta}_1 a_1 - \hat{\eta}_2 a_2 - \left( \mathcal{K}_{\bY}(\eta_1, \tilde{\eta}_2(\eta_1)) -\eta_1 a_1 - \tilde{\eta}_2(\eta_1) a_2 \right),\\
         -\frac{1}{2} \hat{\omega}_2^2 &=& \mathcal{K}_{\bY}(0, \tilde{\eta}_2(0)) - \tilde{\eta}_2(0) a_2, \\
        -\frac{1}{2} (\omega_2 - \hat{\omega}_2)^2 &=& \mathcal{K}_{\bY}(\eta_1, \tilde{\eta}_2(\eta_1)) -\eta_1 a_1 - \tilde{\eta}_2(\eta_1) a_2 - \left( \mathcal{K}_{\bY}(\eta_1,\eta_2) -\eta_1 a_1 - \eta_2 a_2\right).
    \end{eqnarray*}
 The sign of $\bomega$ is chosen for $\omega_i$ to be increasing in $\eta_i$. By the inverse function theorem, there exists an inverse function $\boldeta(\bomega)$.
 To identify the pole after a change of variable, define a function $\tilde{\omega}_2(\omega_1)$ to be the value of $\omega_2$ that makes $\eta_2$ zero when $\omega_1$ is fixed, that is,
    \begin{equation*}
        \eta_2(\omega_1,\tilde{\omega}_2(\omega_1)) = 0.
    \end{equation*}
 Since $\omega_1$ is defined not to depend on $\eta_2$, the determinant of $|\partial \bomega / \partial \boldeta |$ is the product of its diagonals. We can now rewrite (\ref{e-inv2}) as
    \begin{eqnarray}\label{e-inv3}
        \mathsf{E}\left[\overline{X} {\bf 1}_{[\overline{\bY} \geq \ba]}\right] &=& \int_{\hat{\bomega}-i\infty}^{\hat{\bomega}+i\infty} \frac{\exp \left[ n q(\omega_1, \omega_2)\right]}{(2\pi i)^2 \ \eta_1 \eta_2} \ \mathcal{K}_\gamma(\eta_1(\omega_1), \eta_2(\omega_1, \omega_2)) \cdot \frac{\partial \eta_1}{\partial \omega_1}\frac{\partial \eta_2}{\partial \omega_2} \ \rd \bomega  \nonumber \\
        &=& \int_{\hat{\bomega}-i\infty}^{\hat{\bomega}+i\infty} \frac{\exp \left[ n  q(\omega_1, \omega_2) \right]}{(2\pi i)^2 \ \omega_1 (\omega_2 - \tilde{\omega}_2(\omega_1))} \ \mathcal{K}_\gamma(\eta_1, \eta_2) \  F(\eta_1, \eta_2) \ \rd \bomega
    \end{eqnarray}
 where
    \begin{equation*}
        q(\omega_1, \omega_2) = \frac{1}{2}\omega_1^2 +  \frac{1}{2}\omega_2^2 - \hat{\omega}_1 \omega_1 - \hat{\omega}_2 \omega_2
    \end{equation*}
 and
    \begin{equation*}\label{def:f1}
        F(\eta_1, \eta_2) =  \frac{\omega_1}{\eta_1}\frac{\partial \eta_1}{\partial \omega_1} \cdot \frac{\omega_2 - \tilde{\omega}_2(\omega_1)}{\eta_2}\frac{\partial \eta_2}{\partial \omega_2}.
    \end{equation*}

 We closely follow the program set by \cite{Kolassa Li:10} and \cite{Li:08}, but we face additional difficulties because of the term $\mathcal{K}_\gamma(\eta_1, \eta_2)$.
 Decompose $F(\eta_1, \eta_2)$ as
    \begin{equation*}
        F = H^0 + H^1 + H^2 + H^{12},
    \end{equation*}
 where $H^0 (\eta_1, \eta_2) = F(0,0), \ H^1(\eta_1, \eta_2) = F(\eta_1,0) - F(0,0), \ H^2(\eta_1, \eta_2) = F(0,\eta_2) - F(0,0),$ and $H^{12}(\eta_1, \eta_2) = F(\eta_1, \eta_2) - F(\eta_1, 0) - F(0,\eta_2) + F(0,0)$.
 It is proved that $F(0,0) = 1$ and that
    \begin{equation*}
        H^0, \ \frac{H^1}{\omega_1}, \ \frac{H^2}{\omega_2 - \tilde{\omega}_2(\omega_1)}, \mbox{ and } \frac{H^{12}}{\omega_1(\omega_2 - \tilde{\omega}_2(\omega_1))}
    \end{equation*}
 are analytic.
 Then, (\ref{e-inv3}) is decomposed into four terms denoted by $I^0, I^1, I^2$, and $I^{12}$, depending on the respective superscript of $H$.

 In order to compute each integral, we impose the assumption that $\mathcal{K}_\gamma(\eta_1, \eta_2)$ is analytic in a neighborhood of $(\hat{\eta}_1, \hat{\eta}_2)$ containing $(\hat{\eta}_1, 0)$, $(0, \tilde{\eta}_2(0))$, and $(0,0)$.
 The simplest part, $I^{12}$ and $I^{2}$, can be obtained by applying multivariate Watson's lemma after modifying the integrand of $I^2$.
 High-order terms of (\ref{lem:I2+I12}) can be computed, but since the order of $I^0$ and $I^1$ is limited to $O\left(n^{-1}\right)$, we present the result up to $O\left(n^{-1}\right)$.

 \begin{lem}\label{lem:I2}
 The sum of integrals $I^1$ and $I^{12}$ is expanded as
    \begin{eqnarray}\label{lem:I2+I12}
        I^2 + I^{12} &=& \frac{1}{\sqrt{n
        }} \bar{\Phi}(\sqrt{n}\hat{\omega}_1)\phi(\sqrt{n}\hat{\omega}_2)\mathcal{K}_\gamma(0, \tilde{\eta}_2(0)) \\
        && \times \left[\frac{1}{\tilde{\eta}_2(0) \sqrt{ \mathcal{K}_{\bY}^{22}(0,\tilde{\eta}_2(0))}} - \frac{1}{\hat{\omega}_2} \right] +  O\left(n^{-1}\right). \nonumber
    \end{eqnarray}
 \end{lem}
 \proof{ See Appendix \ref{app:I2}.
 \endproof}

 For $I^0$ and $I^1$, we do a change of variable with $(v_1, v_2) = (\omega_1, \omega_2 - \tilde{\omega}_2(\omega_1))$ and set $\hat{v}_1 = \hat{\omega}_1$, $\hat{v}_2 = \hat{\omega}_2 - \tilde{\omega}_2(\hat{\omega}_1)$, and $\tilde{v}_2(0) = \hat{w}_2$.
 Let $\tilde{\mathcal{K}}_\gamma(v_1, v_2)$ denote the function $\mathcal{K}_\gamma$ in terms of $(v_1, v_2)$.
 After the change of variable $(\omega_1, \omega_2 - \tilde\omega_2(\omega_1)) \mapsto (v_1, v_2)$, $I^0$ and $I^1$ are written as
 \begin{eqnarray}
   I^0 &=& \int_{\hat \bv - i\infty}^{\hat\bv + i\infty} \frac{\exp[ng(v_1,v_2)]}{(2\pi i)^2} \frac{1}{v_1v_2} \tilde\cK_\gamma(v_1,v_2)\rd \bv \label{I0} \ \mbox{ and }\\
   I^1 &=& \int_{\hat\bv - i\infty}^{\hat\bv + i \infty} \frac{\exp[ng(v_1,v_2)]}{(2\pi i)^2} \frac{1}{v_2} \tilde\cK_\gamma(v_1,v_2) h(v_1) \rd \bv, \label{I1}
 \end{eqnarray}
 respectively, where $\bv = (v_1, v_2)$, $g(v_1, v_2) = v_1^2 /2 + (v_2 + \tilde{\omega}_2(v_1))^2 / 2 - \hat{\omega}_1 v_1 - \hat{\omega}_2 (v_2 + \tilde{\omega}_2(v_1))$, and an analytic function
    \begin{equation*}
        h(v_1) = \frac{F(\eta_1(v_1),0) - 1}{v_1} = \frac{1}{\eta_1(v_1)} \frac{\rd\eta_1}{\rd v_1} - \frac{1}{v_1}.
    \end{equation*}

 Now, we decompose $\tilde \cK_\gamma(v_1,v_2)$ into four terms as
 \begin{eqnarray}\label{Ktilde4}
    \tilde\cK_\gamma(v_1,v_2) &=& \tilde \cK_\gamma(0,0)  + \Big[\tilde \cK_\gamma(v_1,0) - \tilde \cK_\gamma(0,0) \Big] + \Big[\tilde \cK_\gamma(0,v_2) - \tilde \cK_\gamma(0,0) \Big] \nonumber \\
    && + \Big[\tilde \cK_\gamma(v_1,v_2) - \tilde \cK_\gamma(v_1,0) - \tilde \cK_\gamma(0,v_2) + \tilde \cK_\gamma(0,0)\Big].
 \end{eqnarray}
 By the assumption on $\cK_\gamma(\eta_1, \eta_2)$ and by the composition theorem of complex variables, there exists a region $\cA$ such that $\tilde{\mathcal{K}}_\gamma(v_1, v_2)$ is analytic in $\cA$ and $\cA$ contains $(\hat{v}_1, \hat{v}_2)$, $(\hat{v}_1, 0)$, $(0, \tilde{v}_2(0))$, and $(0,0)$.
 Partial derivatives $\partial \tilde\cK_\gamma(v_1,v_2) / \partial v_1$, $\partial \tilde\cK_\gamma(v_1,v_2) / \partial v_2$, and $\partial^2 \tilde\cK_\gamma(v_1,v_2) / \partial v_1 \partial v_2$ are also analytic in $\cA$.
 Furthermore,
 $$ \frac{\tilde \cK_\gamma(v_1,0) - \tilde\cK_\gamma(0,0)}{v_1}, \frac{\tilde \cK_\gamma(0,v_2) - \tilde \cK_\gamma(0,0)}{v_2},$$
 and
 $$\frac{\tilde \cK_\gamma(v_1,v_2) - \tilde \cK_\gamma(v_1,0) - \tilde \cK_\gamma(0,v_2) + \tilde \cK_\gamma(0,0)}{v_1 v_2}$$
 are analytic as well.
 By plugging (\ref{Ktilde4}) into (\ref{I0}) and (\ref{I1}), the integral $I^0 + I^1$ can be rewritten as
    \begin{eqnarray}\label{I01-d}
      I^0 + I^1 &=& \tilde{\mathcal{K}}_\gamma(0,0) \int_{\hat{\mathbf{v}}-i\infty}^{\hat{\mathbf{v}}+i\infty} \frac{\exp \left[ n g(v_1, v_2)\right]}{(2 \pi i)^2} \frac{1}{v_1 v_2} \rd \mathbf{v}
      + \int_{\hat{\mathbf{v}}-i\infty}^{\hat{\mathbf{v}}+i\infty} \frac{\exp \left[ n g(v_1, v_2)\right]}{(2 \pi i)^2} \frac{1}{v_2} k_1(v_1) \rd \mathbf{v} \nonumber \\
            && + \int_{\hat{\mathbf{v}}-i\infty}^{\hat{\mathbf{v}}+i\infty} \frac{\exp \left[ n g(v_1, v_2)\right]}{(2 \pi i)^2} \frac{1}{v_1} k_2(v_2) \rd \mathbf{v} + O\left(n^{-1}\right)
    \end{eqnarray}
    where
    $$ k_1(v_1) = \frac{\tilde{\mathcal{K}}_\gamma(v_1,0) - \tilde{\mathcal{K}}_\gamma(0,0)}{v_1} +  \tilde{\mathcal{K}}_\gamma(v_1,0) \left( \frac{1}{\eta_1(v_1)} \frac{\rd\eta_1}{\rd v_1} - \frac{1}{v_1} \right)$$
    and
    $$ k_2(v_2) = \frac{\tilde{\mathcal{K}}_\gamma(0, v_2) - \tilde{\mathcal{K}}_\gamma(0,0)}{v_2}. $$
 The terms with analytic integrands disappear as we apply multivariate Watson's lemma since their contributions are of order $O\left(n^{-1}\right)$.

 The importance of the decomposition (\ref{I01-d}) lies in that $I^0$ and $I^1$ are now the sum of certain integrals such that each term can be treated separately via, e.g., the method proposed in \cite{Kolassa:03}. The special case for a bivariate random vector is well described in Chapters 3 and 5 of \cite{Li:08}. To approximate the first term in \eqref{I01-d}, the author approximates $\tilde{\omega}_2(\omega_1) / \omega_1$ by a linear function of $\omega_1$, namely $\tilde{\omega}_2(\omega_1) / \omega_1 = b_0 + b_1(\omega_1 - \hat{\omega}_1)$ since $\tilde{\omega}_2(\omega_1) $ is usually intractable. Then it is proved that the derived saddlepoint expansion using the linear function is equivalent to the saddlepoint expansion without the linear approximation up to the order $O\left(n^{-1}\right)$. As for the second and third integrals, one can expand $g(v_1, v_2)$ about $(\hat{v}_1, \hat{v}_2)$ and integrate termwise, dropping the terms that contribute the error of $O(n^{-r})$ with $r>1$. The same treatments applied to $I^{\{1\}}$ and $I^{\{2\}}$ in \cite{Li:08} lead us to saddlepoint expansions of the second and third integrals, respectively. We do not report the procedure in detail, but summarize the outcome below.

 In the rest of this section, we define some auxiliary variables that appear in our expansion.
 Let $\check{\omega}_2 = \tilde{\omega}_2(\hat{\omega}_1)$ and let $\check{\omega}'_2$ and $\check{\omega}''_2$ be the first and second derivative of $\tilde{\omega}_2$ evaluated at $\hat{\omega}_1$. They can be specifically computed as
 \begin{eqnarray*}
    \check{\omega}_2 &=& \hat{\omega}_2 + {\rm sign}(-\hat{\eta}_2) \sqrt{-2 \left[ \cK_{\bY}(\hat{\eta}_1, \hat{\eta}_2) -  \hat{\eta}_1 a_1  - \hat{\eta}_2 a_2 - ( \cK_{\bY}(\hat{\eta}_1, 0) -  \hat{\eta}_1 a_1  ) \right]},\\
    \check{\omega}'_2 &=& \left. \left( \mathcal{K}_{\bY}^{1}(\hat{\eta}_1,0) - a_1 \right) \left.\frac{\rd\eta_1}{\rd\omega_1}\right|_{\hat{\omega}_1} \right/ (\check{\omega}_2 - \hat{\omega}_2), \mbox{ and }
 \end{eqnarray*}
 \begin{eqnarray*}
    \check{\omega}''_2 &=& \left[ \Big(\mathcal{K}_{\bY}^{11}(\hat{\eta}_1,0)-\mathcal{K}_{\bY}^{11}(\hat{\eta}_1,\hat{\eta}_2) -\mathcal{K}_{\bY}^{12}(\hat{\eta}_1,\hat{\eta}_2)\tilde{\eta}'_2(\hat{\eta}_1)\Big)\left( \left.\frac{\rd\eta_1}{\rd\omega_1}\right|_{\hat{\omega}_1} \right)^2 \right. \\
    && \left.\left. + (\mathcal{K}_{\bY}^{1}(\hat{\eta}_1,0)-a_1) \left.\frac{\rd^2\eta_1}{\rd\omega_1^2}\right|_{\hat{\omega}_1}-(\check{\omega}'_2)^2 \right] \right/ (\check{\omega}_2 - \hat{\omega}_2).
 \end{eqnarray*}
 Here,
 \begin{eqnarray*}
    \left.\frac{\rd\eta_1}{\rd\omega_1}\right|_{\hat{\omega}_1}  &=& \sqrt{\frac{1}{\mathcal{K}_{\bY}^{11}(\hat{\eta}_1,\hat{\eta}_2) + \mathcal{K}_{\bY}^{12}(\hat{\eta}_1,\hat{\eta}_2)\tilde{\eta}'_2(\hat{\eta}_1) }}
 \mbox{ with }
 \tilde{\eta}'_2(\hat{\eta}_1) = - \frac{\mathcal{K}_{\bY}^{12}(\hat{\eta}_1,\hat{\eta}_2)}{\mathcal{K}_{\bY}^{22}(\hat{\eta}_1,\hat{\eta}_2)}; \mbox{ and } \\
    \left.\frac{\rd^2\eta_1}{\rd\omega_1^2}\right|_{\hat{\omega}_1} &=& -  \Big[ \Big(\mathcal{K}_{\bY}^{111}(\hat{\eta}_1,\hat{\eta}_2) + 2\mathcal{K}_{\bY}^{112}(\hat{\eta}_1,\hat{\eta}_2)\tilde{\eta}'_2(\hat{\eta}_1) + \mathcal{K}_{\bY}^{122}(\hat{\eta}_1,\hat{\eta}_2)\tilde{\eta}'_2(\hat{\eta}_1)^2 \\
    && + \left. \left. \left. \mathcal{K}_{\bY}^{12} (\hat{\eta}_1,\hat{\eta}_2)\tilde{\eta}''_2(\hat{\eta}_1)\right) \left( \left.\frac{\rd\eta_1}{\rd\omega_1}\right|_{\hat{\omega}_1} \right)^2 \right] \right/\Big(3(\mathcal{K}_{\bY}^{11}(\hat{\eta}_1,\hat{\eta}_2) +\mathcal{K}_{\bY}^{12}(\hat{\eta}_1,\hat{\eta}_2)\tilde{\eta}'_2(\hat{\eta}_1))\Big)
 \end{eqnarray*}
 with
 $\tilde{\eta}''_2(\hat{\eta}_1) = -\Big[\mathcal{K}_{\bY}^{112}(\hat{\eta}_1,\hat{\eta}_2) + 2\mathcal{K}_{\bY}^{122}(\hat{\eta}_1,\hat{\eta}_2)\tilde{\eta}'_2(\hat{\eta}_1) + \mathcal{K}_{\bY}^{222}(\hat{\eta}_1,\hat{\eta}_2)\tilde{\eta}'_2(\hat{\eta}_1)^2 \Big] / \mathcal{K}_{\bY}^{22}(\hat{\eta}_1,\hat{\eta}_2)$.
 Then we have $b_0 = \check{\omega}'_2 - \check{\omega}''_2 ~\hat{\omega}_1/2$ and $b_1 = \check{\omega}''_2/2$.
 Moreover,
 let $\hat{x}=\sqrt{n}(\hat{\omega}_1+b_0 \hat{\omega}_2)/\sqrt{1+b_0^2}$, $\hat{y}=\sqrt{n}~\hat{\omega_2}$, $\hat{\rho}=b_0/\sqrt{1+b_0^2}$, $\hat{t} = \sqrt{n}\sqrt{1+b_0^2} \hat{\omega}_1$,
 and  $\check{g}=(\check{\omega}_2-\check{\omega}_2^\prime\hat{\omega}_1)\left(\check{\omega}_2/2 - \check{\omega}_2^\prime\hat{\omega}_1/2-\hat{\omega}_2 \right)$.

 The extension of Theorem \ref{thm2} for $d=2$ is presented by summarizing the above arguments in Theorem \ref{thm3}.
 \begin{thm}\label{thm3}
    Let $\hat{\boldeta}$ solve the equation (\ref{spaeqn-multi}) with $\hat\eta_i > 0$ for $i=1,2$, and suppose that $\mathcal{K}_\gamma(\boldeta)$ is analytic in a neighborhood of $\hat{\boldeta}$ containing $(\hat{\eta}_1, 0)$, $(0, \tilde{\eta}_2(0))$, and $(0,0)$.
    With all the notation defined above,
    $\mathsf{E}\left[\overline{X} {\bf 1}_{[\overline{\bY} \geq \ba]}\right] $ of a continuous random vector $(X,\mathbf{Y}) \in \mathbb{R}^3$ can be approximated via saddlepoint techniques by
    \begin{eqnarray*}
        \mathsf{E}\left[\overline{X} {\bf 1}_{[\overline{\bY} \geq \ba]}\right]  &=& \mathsf{E}[X] \bar{\Phi}(\hat{x}, \hat{y}, \hat{\rho}) + \frac{1}{\sqrt{n}} \left\{ \frac{\mathsf{E}[X] b_1}{1+b_0^2}\phi(\hat{x}) \left[ \sqrt{1-\hat{\rho}^2}\left(\hat{x} - \hat{t} \right)~\phi\left( \frac{\hat{y}-\hat{\rho}\hat{x}}{\sqrt{1-\hat{\rho}^2}} \right)  \right. \right. \\
        && - \left. \left(\hat{\rho}+\hat{x}\hat{y}-\hat{\rho}\hat{x}^2 - \hat{y}\hat{t} +\hat{\rho}\hat{x}\hat{t}\right) \bar{\Phi}\left( \frac{\hat{y}-\hat{\rho}\hat{x}}{\sqrt{1-\hat{\rho}^2}} \right) \right] + \mathcal{K}_\gamma(0, \tilde{\eta}_2(0)) \\
        && \times \left[\frac{1}{\tilde{\eta}_2(0) \sqrt{ \mathcal{K}_{\bY}^{22}(0,\tilde{\eta}_2(0))}} - \frac{1}{\hat{\omega}_2} \right] \phi(\sqrt{n}\hat{\omega}_2) \bar{\Phi}(\sqrt{n}\hat{\omega}_1) + \exp\left[n\check{g}\right] \\
        && \times \left[ \frac{k_1\left(\hat{\omega}_1\right)}{\sqrt{1+(\check{\omega}_2^\prime)^2}} \phi\left(\frac{\sqrt{n}[(1+(\check{\omega}_2^\prime)^2)\hat{\omega}_1 + \check{\omega}_2^\prime(\hat{\omega}_2-\check{\omega}_2)]}{\sqrt{1+(\check{\omega}_2^\prime)^2}}\right)\bar{\Phi}\left(\frac{\sqrt{n}(\hat{\omega}_2-\check{\omega}_2)}{\sqrt{1+(\check{\omega}_2^\prime)^2}}\right)   \right. \\
        && + k_2\left(\hat{\omega}_2\right) \phi \left(\sqrt{n} (\check{\omega}_2' \hat{\omega}_1 + \hat{\omega}_2 - \check{\omega}_2 )\right) \bar{\Phi}\left( \sqrt{n}\hat{\omega}_1\right)  \Big]\Big\} + O\left(n^{-1}\right)
    \end{eqnarray*}
    where
    \begin{equation*}
        k_1\left(\hat{\omega}_1\right) = \frac{\mathcal{K}_\gamma(\hat{\eta}_1,0) - \mathcal{K}_\gamma(0,0)}{\hat{\omega}_1} + \mathcal{K}_\gamma(\hat{\eta}_1,0) \left[ \frac{1}{\hat{\eta}_1 \sqrt{ \mathcal{K}_{\bY}^{11}(\hat{\eta}_1,\hat{\eta}_2) + \mathcal{K}_{\bY}^{12}(\hat{\eta}_1,\hat{\eta}_2)\tilde{\eta}'_2(\hat{\eta}_1) }} - \frac{1}{\hat{\omega}_1} \right]
    \end{equation*}
    and
    \begin{equation*}
        k_2\left(\hat{\omega}_2\right) = \frac{\mathcal{K}_\gamma(0,\tilde{\eta}_2(0)) - \mathcal{K}_\gamma(0,0)}{\hat{\omega}_2}
    \end{equation*}
    for $\hat{\boldeta} > \mathbf{0}$.
    Here, $\bar{\Phi}({x}, {y}, {\rho}) = 1 - \Phi({x}, {y}, {\rho})$ with the CDF $\Phi({x}, {y}, {\rho})$ of a bivariate standard normal variable $\mathcal{N}(0, 0, 1, 1, \rho)$.
 \end{thm}

 \begin{rem}
 We omit the cases where $\hat{\boldeta} = 0$ or at least one component of $\hat{\boldeta}$ is negative due to complexity. However, both can be argued just as bivariate saddlepoint approximations, after applying the decomposition (\ref{I01-d}).
 \end{rem}
 
  \section{Applications in risk management}\label{ch5}

 Saddlepoint techniques have been successfully applied in various problems of quantitative finance such as vanilla option pricing, portfolio risk measurements. The newly developed saddlepoint approximations allow us to extend the applicability to other important problems in risk management. In particular, we consider fast and accurate computations of risk and option sensitivities which are indispensable in responsive decision making.

 First of all, we consider a random portfolio loss $L$, and compute the sensitivities of certain risk metrics utilizing Theorems \ref{thm1} and \ref{thm2}.
 We particularly investigate Euler contributions to risk measures in Section \ref{sec:ch5-app1} and risk sensitivities with respect to an input parameter under a delta-gamma portfolio model in Section \ref{sec:ch5-app2}.
 The second application is on option sensitivities. This exercise is done under two different asset pricing models as described in Section \ref{sec:ch5-app3}.
 Numerical illustrations shall confirm the accuracy and effectiveness of saddlepoint approximations.

 \subsection{VaR and CVaR risk contribution}\label{sec:ch5-app1}

 Suppose that there is a portfolio with continuous random loss $L$, consisting of $m$ assets or sub-portfolios $L_i$'s with $u_i$ units of asset (portfolio) $i$ for $i = 1, \ldots, m$, so that $L = \sum_{i=1}^m u_i L_i$. For a risk measure, say $\nu(L)$, it is important to know how much the sub-portfolio $L_i$ contributes to $\nu(L)$ from a risk management point of view.
 Risk measures of our interest are the most frequently used measures, namely, value-at-risk (VaR) $v_\alpha$, a quantile function of the distribution of $L$, and conditional value-at-risk (CVaR) $c_\alpha$, also called expected shortfall (ES). Fix $\alpha \in (0,1)$, typically taken to be $0.95$ or $0.99$. Then $v_\alpha$ and $c_\alpha$ are given by $v_\alpha = \inf \{ l | {\sf P}(L \leq l) \geq \alpha \}$ and $c_\alpha =  \mathsf{E}[L | L \geq v_\alpha]$. If necessary, we write $v_\alpha(L)$ or $c_\alpha(L)$ to specify the underlying random loss variable.

 For a risk measure that is homogeneous of degree 1 and differentiable in an appropriate sense, the Euler allocation principle can be applied. We refer the reader to \cite{Tasche:08} for more information where the author defines the Euler contributions to VaR and CVaR as
 $v_\alpha(L_i | L) = \mathsf{E}[L_i | L = v_\alpha]$ and
 $c_\alpha(L_i | L) = \mathsf{E}[L_i | L \geq v_\alpha]$.

 As such risk metrics have drawn much attention from researchers and practitioners, saddlepoint approximations to VaR and CVaR risk contributions have been studied in the literature. For example, see \cite{Martin et al:01} or \cite{Muromachi:04}. The VaR risk contribution formula in \cite{Martin et al:01} is simple to apply and is nothing but the first order approximation. On the other hand, the approximations provided in \cite{Muromachi:04} are rather complex to compute. In particular, the expansions make use of an auxiliary function which acts like a CGF, and thus it is difficult to guarantee the existence of saddlepoints.

 \subsubsection{A portfolio composed of correlated normals}

 Suppose that random losses $\{ L_i \}_{i=1,\cdots,m}$ follow a multivariate normal distribution $\mathcal{N}(\bmu, \bSigma)$ with an $m$-dimensional mean vector $\bmu = (\mu_1, \cdots, \mu_m)^\top$ and an $m\times m$ covariance matrix $\bSigma$ whose entries are $\Sigma_{ii} = \sigma_i^2$ and $\Sigma_{ij} = \Sigma_{ji} = \rho_{ij}\sigma_i \sigma_j$ with $\rho_{ij} = \rho_{ji}$.
 We apply Theorems \ref{thm1} and \ref{thm2} to the Euler contributions with $n=1$. The resulting formulas are actually the same as the true values:
 \begin{eqnarray*}
    \mathsf{E}[L_i | L = v_\alpha ] &=& \mu_i + \frac{\bu^\top \bSigma^{i}}{ \bu^\top \bSigma \bu}(v_\alpha - \bu^\top \bmu ), \;\; \mbox{ and }\\
    \mathsf{E}[L_i | L \geq v_\alpha ] &=& \mu_i + \frac{\bu^\top \bSigma^{i}}{ \sqrt{\bu^\top \bSigma \bu}} \cdot \frac{\phi(\hat{\omega})}{1- \Phi(\hat{\omega})}
 \end{eqnarray*}
 where $\bSigma^{i}$ is the $i$-th column of $\bSigma$, $\bu =(u_1, \cdots, u_m)^\top$ and $\hat{\omega} = (v_\alpha - \bu^\top \bmu )/\sqrt{\bu^\top \bSigma \bu}$.

 For comparison, we note that the approach of \cite{Martin et al:01} yields the same result whereas
 Muromachi's formula for the VaR contribution results in
 $$
 \sqrt{\frac{\bu^\top \bSigma \bu}{\mathcal{K}_{M}''(\hat\eta_M)}} \exp \left[ \frac{(\bu^\top \bmu - v_\alpha)^2}{2 \bu^\top \bSigma \bu} + \mathcal{K}_{M}(\hat\eta_M) - \hat\eta_M v_\alpha \right] \left( 1 + \frac{1}{8} \hat{\rho}_{4,M} - \frac{5}{24} \hat{\rho}_{3,M}^2 \right).
 $$
 Here, $\cK_M$ is defined as $\cK_L + \log(\partial \cK_L(\eta)/\partial u_i) - \log \eta$, different from the CGF $\cK_L$ of $L$. In this example, it is given by
 $$
 \mathcal{K}_M (\eta)= \bu^\top \bmu \eta + \frac{1}{2}\bu^\top \bSigma \bu\eta^2 + \log \left[ \mu_i \eta + \left(\sum_{k \neq i} u_k\rho_{ik}\sigma_k\sigma_i + u_i\sigma_i^2\right) \eta^2 \right] - \log \eta.
 $$
 Moreover,
 $\hat{\eta}_M$ is the saddlepoint of $\mathcal{K}_M$, that is, the solution of the following cubic polynomial equation
 \begin{align*}
    & \Big(\bu^\top \bSigma \bu + \bu^\top \bmu- v_\alpha\Big) \bu^\top \bSigma^i  \eta^3 + \Big(\mu_i \left(\bu^\top \bSigma \bu + \bu^\top \bmu- v_\alpha \right) + 2\bu^\top \bSigma^i\Big) \eta^2 \\
    & + \left(\mu_i + \bu^\top \bSigma^i\right)\eta + \mu_i = 0.
 \end{align*}
 Lastly, $\hat\rho_{r,M}$ is the standardized cumulant for $\cK_M$ of order $r$ evaluated at $\hat\eta_M$.

 \subsubsection{A portfolio of proper generalized hyperbolic distributions}

 Consider a proper generalized hyperbolic (GH) distribution which is a GH distribution with the restricted range of parameters $\lambda \in \mathbb{R}, \alpha > 0, \beta \in (-\alpha, \alpha), \delta > 0,$ and $\mu \in \mathbb{R}$. This excludes some cases such as variance gamma distribution, but still continues to nest hyperbolic and normal inverse gaussian distributions.
 Let $X \sim \textrm{pGH}(\lambda, \alpha, \beta, \delta, \mu)$ denote a random variable that has a proper GH distribution with the parameter set $(\lambda, \alpha, \beta, \delta, \mu)$.
 The MGF of the proper GH $X$ is expressed as
 \begin{equation*}
     e^{\mu \gamma} \frac{B_\lambda (\delta \sqrt{\alpha^2 - (\beta + \gamma)^2})}{B_\lambda (\delta \sqrt{\alpha^2 - \beta^2})} \left( \frac{\alpha^2 - \beta^2}{\alpha^2 - (\beta + \gamma)^2}\right)^{\lambda/2}.
 \end{equation*}
 Here, $B_{\lambda}(l)$ is the modified Bessel function of the third kind with index $\lambda$ for $l > 0$.

 Let $L_i \sim \textrm{pGH}(\lambda_i, \alpha_i, \beta_i, \delta_i, \mu_i)$ be independent random variables. The target portfolio loss $L$ is given by $L = \sum_{i=1}^m u_i L_i$ where $u_i > 0$ for each $i$. By the scaling property of GH distributions, it is not difficult to check that $u_i L_i$ has the proper GH distribution with the parameter set $(\lambda_i, \alpha_i/u_i, \beta_i/u_i, u_i \delta_i, u_i \mu_i)$ and that the CGF is given by
 \begin{equation*}
    \cK_{u_iL_i}(\eta) = u_i \mu_i \eta + \log B_{\lambda_i}(w_i Q_i(\eta)) - \log B_{\lambda_i}(w_i) - \lambda_i \log Q_i(\eta),
 \end{equation*}
 where $Q_i(\eta) = \sqrt{1 - (2u_i\beta_i \eta + (u_i\eta)^2)/(\alpha_i^2 - \beta_i^2)}$. Finally, $\cK_L(\eta) = \sum_{i=1}^m \cK_{u_iL_i}(\eta)$. Thanks to the relation
 $$
  -2 B'_\lambda(x) = B_{\lambda-1}(x) +  B_{\lambda+1}(x)
  $$
  for $\lambda \in \mathbb{R}$ and $x \in \mathbb{R}^+$, the first derivative of $\cK_{u_iL_i}$ is seen to be
 \begin{equation*}
    \mathcal{K}_{u_iL_i}'(\eta) = u_i \mu_i + \frac{u_i\beta_i + u_i^2\eta}{Q_i(\eta)(\alpha_i^2 - \beta_i^2)} \left\{ \frac{\varsigma_i}{2} \frac{B_{\lambda_i-1}(\varsigma_i Q_i(\eta)) + B_{\lambda_i+1}(\varsigma_i Q_i(\eta))}{B_{\lambda_i}(\varsigma_i Q_i(\eta))} + \frac{\lambda_i}{Q_i(\eta)} \right\}
 \end{equation*}
 where $\varsigma_i = \delta_i \sqrt{\alpha_i^2 - \beta_i^2}$.  The saddlepoint $\hat{\eta}$ needs to be numerically computed by solving $\mathcal{K}'_L(\eta)  = v_\alpha$. The solution is unique in the convergence interval of the CGF of $L$, $$\Big(\max (-\alpha_i/u_i-\beta_i/u_i), \min (\alpha_i/u_i-\beta_i/u_i)\Big).$$

 The VaR or CVaR risk contribution of the portfolio $L$ for the asset $L_i$, $1 \leq i \leq m$, requires to compute the joint CGF of $(L_i, L)$ in order to apply Theorems \ref{thm1} and \ref{thm2}. The joint CGF can be easily derived as
 \begin{eqnarray*}
    \mathcal{K}_{L_i, L}(\gamma,\eta) &=& \sum_{j=1, j \neq i}^m \mathcal{K}_j(\eta) + \mu_i (u_i \eta + \gamma) + \log B_{\lambda_i}(\varsigma_i \tilde{Q}(\gamma,\eta)) \\
    && - \log B_{\lambda_i}(\varsigma_i) - \lambda_i \log \tilde{Q}(\gamma,\eta),
 \end{eqnarray*}
 where
 $$\tilde{Q}(\gamma,\eta) = \sqrt{1 - \frac{2\beta_k(u_i \eta + \gamma) + (u_i \eta + \gamma)^2}{\alpha_i^2 - \beta_i^2}}.$$
 Then a bit of work shows that
 \begin{eqnarray*}
    \frac{\partial}{\partial \gamma} \cK_{L_i, L}(\gamma,\eta) &=& \mu_i + \frac{\beta_i + u_i \eta + \gamma}{\tilde{Q}(\gamma,\eta)(\alpha_i^2 - \beta_i^2) } \\
    && \times \left\{ \frac{\varsigma_i}{2} \frac{B_{\lambda_i-1}(\varsigma_i \tilde{Q}(\gamma,\eta)) + B_{\lambda_i+1}(\varsigma_i \tilde{Q}(\gamma,\eta))}{B_{\lambda_i}(\varsigma_i \tilde{Q}(\gamma,\eta))} + \frac{\lambda_i}{\tilde{Q}(\gamma,\eta)} \right\},
 \end{eqnarray*}
 which yields $\mathcal{K}_\gamma(\eta) = {\mathcal{K}'_{u_iL_i}(\eta)}/{u_i}$ so that $\mathcal{K}_\gamma(\eta)$ is analytic at $\hat{\eta}$.
 By calculating the cumulants needed for Theorems \ref{thm1} and \ref{thm2}, we can obtain the VaR and CVaR risk contributions analytically except for the saddlepoint $\hat\eta$ which can be efficiently found by any root-finding method.

 For the rest of this subsection, we conduct some numerical experiments with an NIG distribution which is a special case of proper GH distributions. The CGF of $u_iL_i$ is reduced to
 $$
 \mathcal{K}_{u_iL_i}(\eta) = u_i \mu_i \eta + \delta_i\left(\sqrt{\alpha_i^2 - \beta_i^2} - \sqrt{\alpha_i^2 - (\beta_i+u_i\eta)^2}\right).
 $$
 The cumulants of $L$ at $\hat{\eta}$ are easily computed. We also have $\mathsf{E}[L_i] = \mu_i + \delta_i \beta_i / \sqrt{\alpha_i^2 - \beta_i^2}$.

 More specifically, we set $m = 3$, $\bu = (0.2, 0.4, 0.4)$ and $L_1 \sim \textrm{pGH}(-1/2, 2, 0.1, 1.8, 0.2)$,  $L_2 \sim \textrm{pGH}(-1/2, 3, 0.3, 0.5, 0.3)$,  and $L_3 \sim \textrm{pGH}(-1/2, 2.5, -0.2, 1, 0.5)$.
 Figures~\ref{fig:RC-NIG-VaR} and \ref{fig:RC-NIG-ES} show the estimated risk contributions of VaR and CVaR for $L_3$.
 We obtain two estimates of VaR $v_\alpha(L)$ using Monte Carlo simulation and saddlepoint techniques, denoted by `MC-VaR' and `SPA-VaR' respectively.
 Then the VaR contribution is computed first using MC-VaR, denoted by `SPA from MC-VaR', and second using SPA-VaR, denoted by `SPA from SPA-VaR'.
 We also plot the approximate VaR contribution, `Martin-SPA', given in \cite{Martin et al:01}.
 For comparison, we compute Monte Carlo estimates based on infinitesimal perturbation analysis, or simply IPA estimates, developed in \cite{Hong:08} using $2 \cdot 10^7$ random outcomes. The batch size for each VaR contribution estimate is set equal to $2 \cdot 10^4$.

 Figure \ref{fig:RC-NIG-VaR} shows that our approximation formulas give very accurate values and that there is a notable difference between Martin-SPA and the others. For a better comparison, the differences between the SPA based estimates and the IPA estimates are shown in the right panel (ii) of Figure \ref{fig:RC-NIG-VaR}. In the monitoring interval of $\alpha$, those absolute or relative differences stay small. For example, the average relative difference between the IPA estimates and SPA from SPA-VaR is $6.0147\times 10^{-3}$. The fluctuating behavior of the difference between the estimates is due to the strong dependence of the IPA estimator on the batch size.
 \begin{figure}[h]
    \centering
    \renewcommand{\thesubfigure}{(i)}
    \subfigure[][]{\includegraphics[width=0.49\linewidth, height=.36\linewidth]{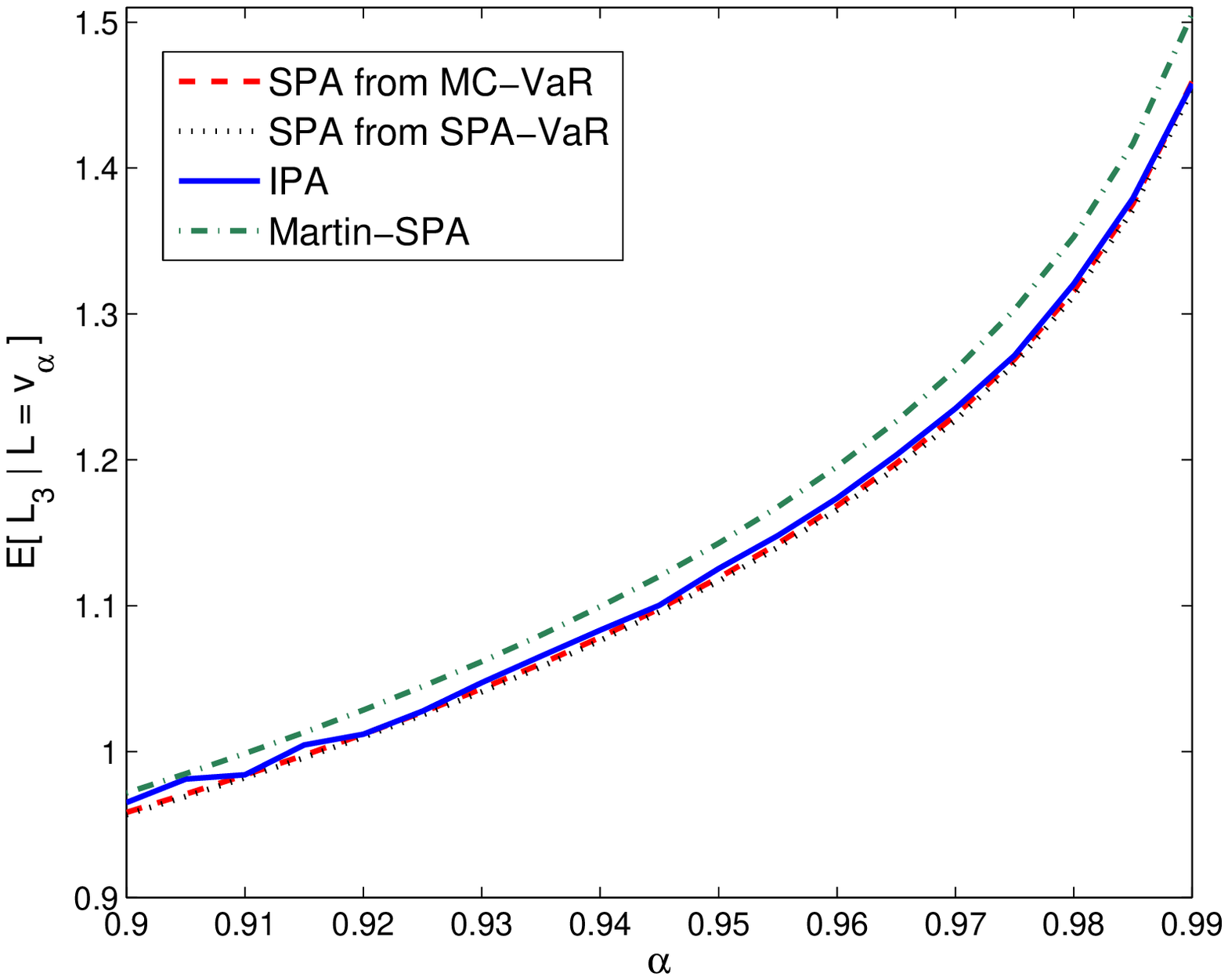}}
    \renewcommand{\thesubfigure}{(ii)}
    \subfigure[][]{\includegraphics[width=0.49\linewidth, height=.36\linewidth]{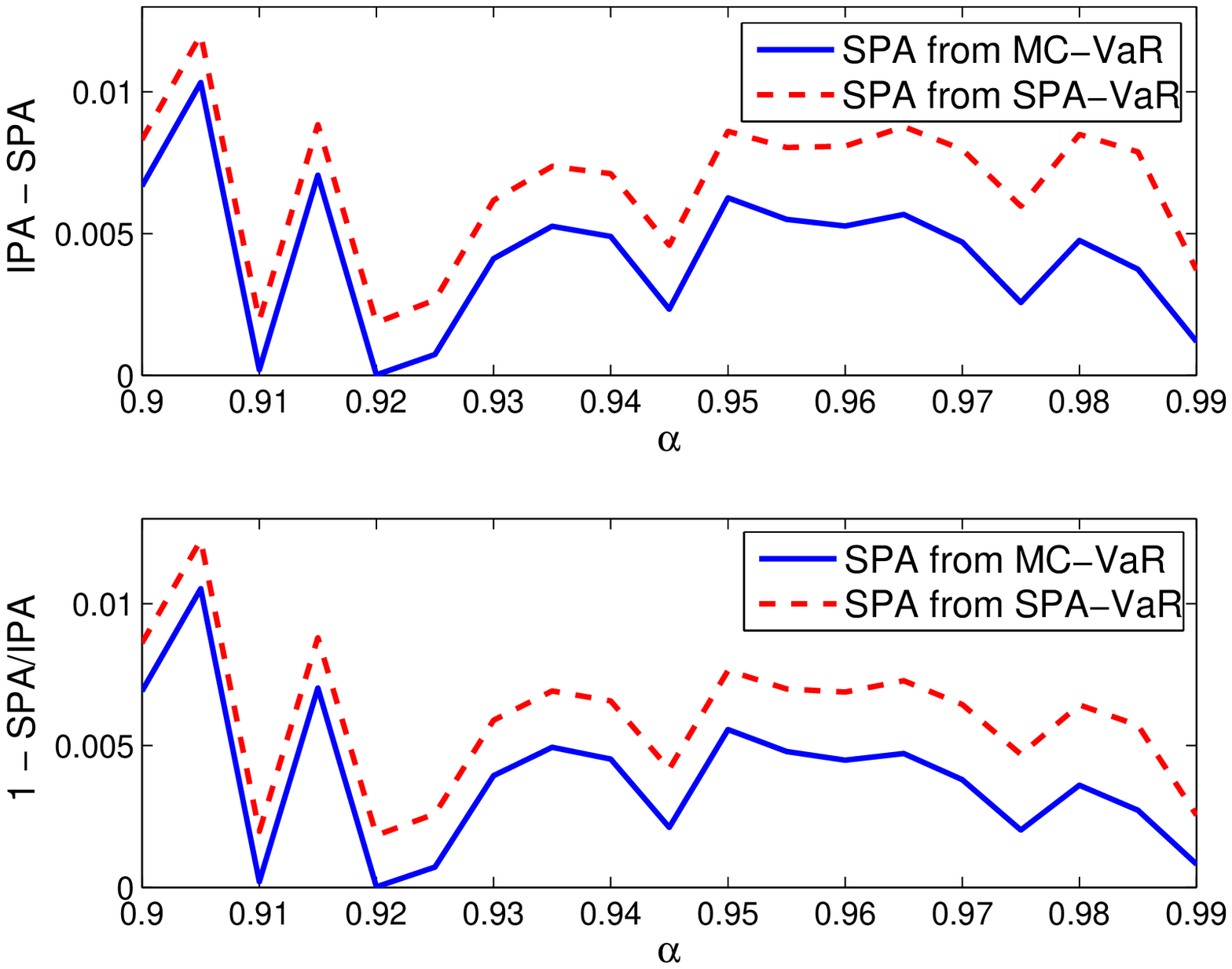}}
    \caption{(i) VaR contribution of $L_3$ over $\alpha$ using SPAs (red, black) and IPA estimator (blue) and (ii) the estimated differences and relative differences of our SPAs to IPA estimator.}\label{fig:RC-NIG-VaR}
    \centering
 \end{figure}

 Figure \ref{fig:RC-NIG-ES} plots the CVaR sensitivities computed by saddlepoint approximations using two VaR estimates, MC-VaR and SPA-VaR, and the results are again denoted by (red) SPA from MC-VaR and (black) SPA from SPA-VaR, respectively.
 As seen from the figure, our SPA formulas from both MC-VaR and SPA-VaR provide highly accurate approximations to the CVaR contribution. For instance, the average relative difference between the IPA estimates and SPA from SPA-VaR is $2.4469\times 10^{-3}$.
 \begin{figure}[h]
    \centering
    \renewcommand{\thesubfigure}{(i)}
    \subfigure[][]{\includegraphics[width=0.49\linewidth, height=.36\linewidth]{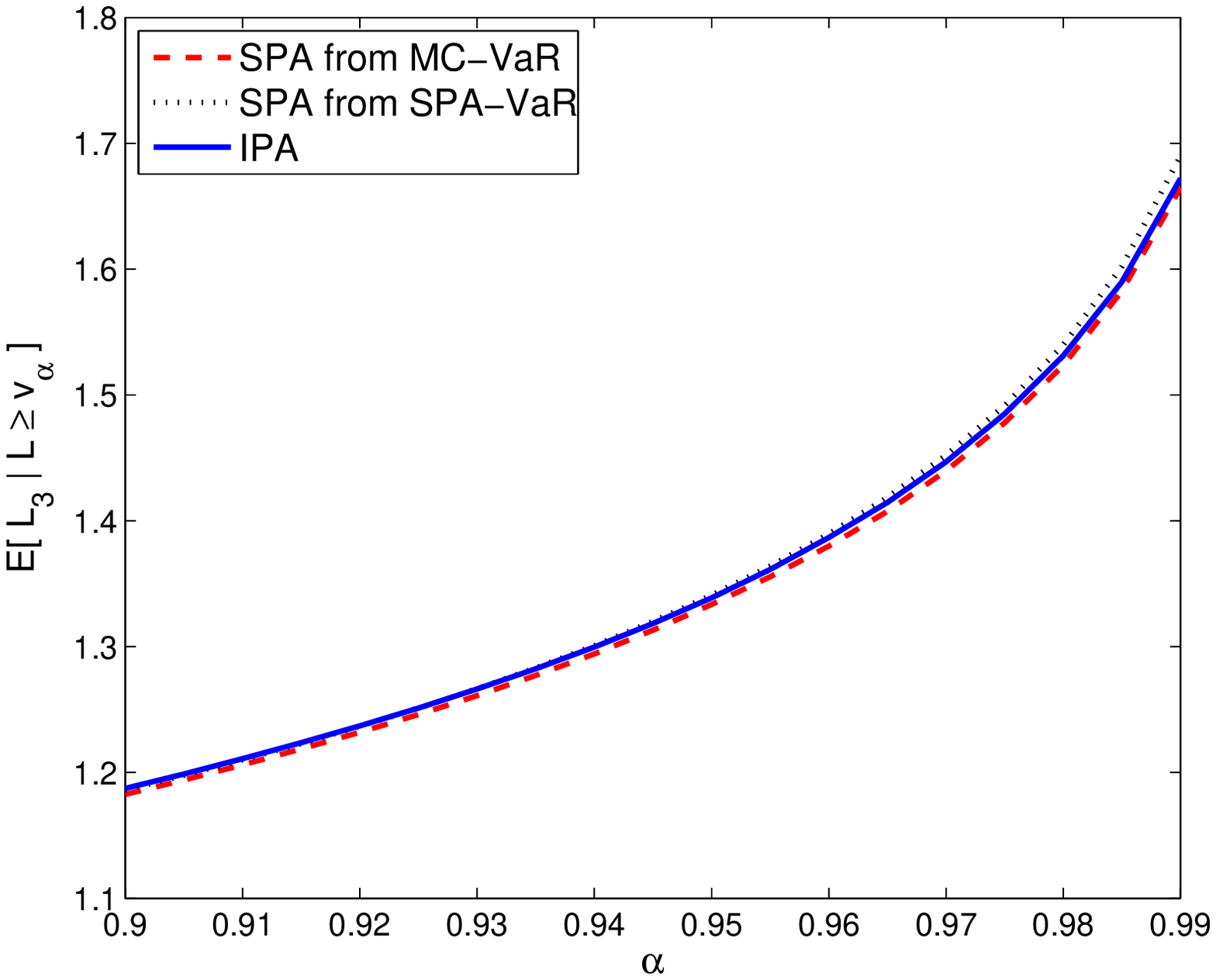}}
    \renewcommand{\thesubfigure}{(ii)}
    \subfigure[][]{\includegraphics[width=0.49\linewidth, height=.36\linewidth]{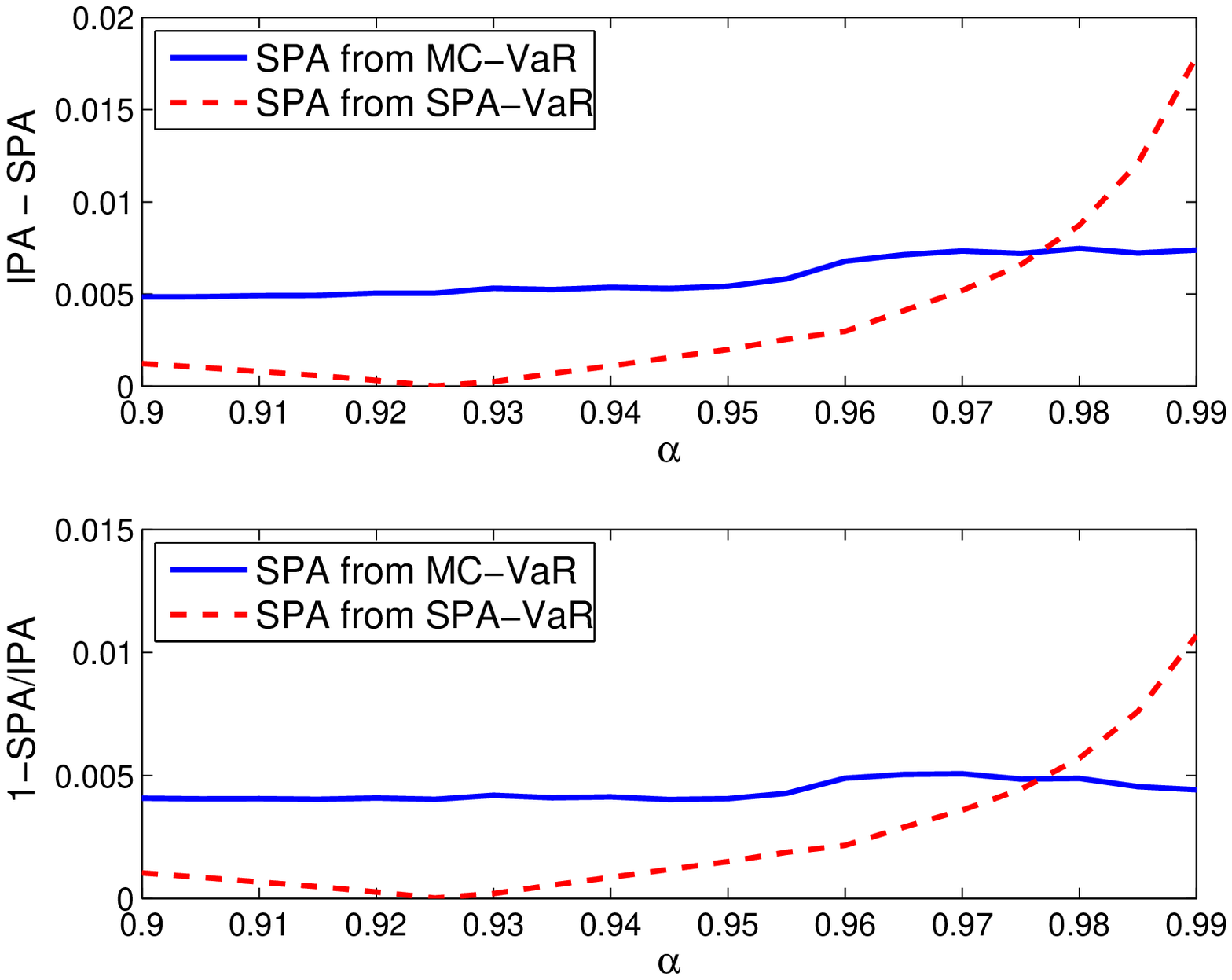}}
    \caption{(i) CVaR contribution of $L_3$ over $\alpha$ using SPAs (red, black) and IPA estimator (blue) and (ii) the estimated differences and relative differences of our SPAs to IPA estimator.}\label{fig:RC-NIG-ES}
    \centering
 \end{figure}

 \subsection{VaR and CVaR sensitivities of delta-gamma portfolios}\label{sec:ch5-app2}

 A delta-gamma portfolio can be understood as a quadratic approximation to portfolio returns and it has been widely employed in quantitative risk management. For example, it is useful in computing VaR of a portfolio loss that could occur in a short period of time. In this section, we extend the existing results on delta-gamma portfolios by computing VaR and CVaR sensitivities with respect to an input parameter.

 \cite{Hong:08} and \cite{Hong Liu:09} show that the sensitivities of $v_\alpha$ and $c_\alpha$ with respect to a general input parameter can be described as  conditional expectations.
 Let the random loss of a portfolio $L(\theta) = \psi(\theta, Z)$ be a function of $\theta$ and a random variable $Z$, where $\theta$ is the parameter with respect to which we differentiate.
 Under certain technical assumptions in \cite{Hong:08}, the VaR sensitivity with respect to $\theta$ can be written as
 \begin{equation*}
    \frac{\partial v_\alpha}{\partial \theta} = \mathsf{E}\left[\frac{\partial \psi}{\partial\theta}(\theta, Z) \Big| \psi(\theta,Z) = v_\alpha\right].
 \end{equation*}
 On the other hand \cite{Hong Liu:09}  prove that  the CVaR sensitivity with respect to $\theta$ is simply
 \begin{equation*}
    \frac{\partial c_\alpha}{\partial \theta} = \mathsf{E}\left[\frac{\partial \psi}{\partial \theta}(\theta, Z) \Big| \psi(\theta,Z) \geq v_\alpha\right].
 \end{equation*}
 as long as certain conditions are met. And the authors develop IPA based estimators using Monte Carlo sampling.

 \subsubsection{Delta-gamma portfolios}

 We first present a setting for a delta-gamma portfolio according to \cite{Feu. Wong:00}. Let a random vector $\bX = (X_1, \cdots, X_m)^\top$ represent the $m$ underlying risk factors in a financial market over a given time period. As often done in the literature, we assume that $\bX$ follows a multivariate normal distribution with mean vector $\bmu$ and covariance matrix $\bSigma$. These parameters are assumed to be known, but in practice they need to be estimated from either historical data or market data. We are concerned with a portfolio loss due to the random factor $\bX$, which we simply denote by $f(\bX)$ for some functional $f$.

 Taking the Taylor expansion of $f(\bX)$ at $\bX=\mathbf{0}$ up to the second order yields a delta-gamma portfolio loss $Y$ for the given time horizon as
 \begin{equation}\label{DG_Y}
    Y = f(\bX) = f(\mathbf{0}) + \ba^\top \bX + \bX^\top \bB \bX,
 \end{equation}
 where $\ba$ is an $m\times1$ column vector and $\bB$ is a symmetric $m \times m$ matrix.
 In order to compute the CGF of $Y$, rewrite $Y$ with zero-mean vector multivariate Gaussian $\bX_0$ as
 $$
 Y = f(\mathbf{0}) + \ba^\top (\bmu + \bX_0) + (\bmu + \bX_0)^\top \bB (\bmu + \bX_0) = c + (\ba + 2\bB \bmu)^\top \bX_0 + \bX_0^\top \bB \bX_0,
 $$
 where $c = f({\bf 0}) + \ba^\top \bmu + \bmu^\top \bB \bmu$.
 Let $\bX_0 = \bH \widetilde{\bZ}$ with an $m \times 1$ column vector $\widetilde{\bZ}$ of independent standard normal random variables using an $m \times m$ matrix $\bH$ such that $\bSigma = \bH \bH^\top$. Performing an eigenvalue decomposition gives us $\bH^\top \bB \bH = \bP \bLambda \bP^\top$, where $\bLambda=$ diag$(\lambda_1,\cdots,\lambda_m)$ is the diagonal matrix of eigenvalues, and $\bP$ is an orthonormal matrix whose $i$-th column is the $i$-th eigenvector associated with the $i$-th eigenvalue $\lambda_i$. This decomposition finally allows us to have
 \begin{eqnarray*}
    Y &=& c + (\ba + 2\bB \bmu)^\top \bH \widetilde{\bZ} + \widetilde{\bZ}^\top \bH^\top \bB \bH \widetilde{\bZ} \nonumber\\
    &=& c + (\ba + 2\bB \bmu)^\top \bH \bP \bP^\top \widetilde{\bZ} + \widetilde{\bZ}^\top \bP^\top \bLambda \bP \widetilde{\bZ} \nonumber \\
    &=& c + \bd^\top \bZ + \bZ^\top \bLambda \bZ
 \end{eqnarray*}
 where $\bd = \bP^\top \bH^\top (\ba + 2\bB \bmu)$ and $\bZ = \bP^\top \widetilde{\bZ}$. Note that $\bZ$ consists of independent standard normal entries $Z_i$ for ${i=1,\cdots,m}$.

 Writing $Y  = c + \sum_{i=1}^m (d_i Z_i + \lambda_i Z_i^2)$ where $d_i$ stands for the $i$-th element of $\bd$, we can compute the MGF $\mathcal{M}_Y(\eta)$ and the CGF $\mathcal{K}_Y(\eta)$ of $Y$ as follows:
 \begin{eqnarray}
    \mathcal{M}_Y(\eta) &=& \left( \prod_{i=1}^m (1 - 2 \lambda_i \eta)\right)^{-1/2} \exp \left( c\eta + \frac{1}{2}\sum_{i=1}^m \frac{d_i^2 \eta^2}{1 - 2 \lambda_i \eta}\right), \mbox{ and } \nonumber \\
    \mathcal{K}_Y(\eta) &=& c\eta -\frac{1}{2}\sum_{i=1}^m \log(1 - 2 \lambda_i \eta) + \frac{1}{2}\sum_{i=1}^m \frac{d_i^2 \eta^2}{1 - 2 \lambda_i \eta}. \label{cgf:DG_Y}
 \end{eqnarray}
 Note that both of them are analytic near the origin and we can explicitly obtain the convergence region.
 The saddlepoint $\hat{\eta}$ of $\mathcal{K}_Y(\eta)$ is obtained by solving $\cK'_Y(\eta) = v_\alpha(Y)$, which turns out to be equivalent to solving an $(n+2)$-th order polynomial equation.
 The existence of a unique saddlepoint in a delta-gamma portfolio is always guaranteed.

 \subsubsection{VaR and CVaR sensitivities with respect to the mean vector}

 In this subsection, we obtain more detailed formulas for risk sensitivities by specifying $\theta$ as the mean vector $\bmu$. In addition to the direct implications that risk sensitivities provide, such computations are helpful in assessing the robustness of the estimates of risk measures when the estimation error of $\mu_i$ is not negligible as pointed out by \cite{Hong Liu:09}.

 The variable of our interest is then
 \begin{eqnarray}\label{partial_Y}
    \frac{\partial Y}{\partial \mu_i}  &=& \frac{\partial c}{\partial \mu_i} + \sum_{k=1}^m \left( \frac{\partial d_k}{\partial \mu_i} Z_k + \frac{\partial \lambda_k}{\partial \mu_i}  Z_k^2 \right)  \\
    &=& a_i + 2 \sum_{k=1}^m b_{ik} \mu_k + \sum_{k=1}^m [2\bP^\top \bH^\top \bB]_{ki} Z_k, \nonumber
 \end{eqnarray}
 where $a_i$ is the $i$-th element of $\ba$, $b_{ik}$ is the $(i,k)$-th component of $\bB$, and $[\bM]_{ki}$ represents the $(k,i)$-th component of a matrix $\bM$.
 The joint CGF of a bivariate random vector ($\partial Y /\partial \mu_i, Y)$ is evaluated using the representation (\ref{partial_Y}) as
 \begin{eqnarray*}
    \cK_{\partial_i Y,Y}(\gamma, \eta) &=&  \left( a_i + 2 \sum_{k=1}^m b_{ik} \mu_k \right) \gamma + c\eta -\frac{1}{2}\sum_{k=1}^m \log(1 - 2 \lambda_k \eta) \nonumber \\
    && + \frac{1}{2}\sum_{k=1}^m \frac{([2\bP^\top \bH^\top \bB]_{ki} \gamma + d_k \eta)^2}{1 - 2 \lambda_k \eta}. \label{cgf:DG_dYY}
\end{eqnarray*}
 Here, we denote $\partial Y/\partial \mu_i$ as $\partial_i Y$ for brevity.
 Furthermore, we directly get
 \begin{eqnarray*}\label{cgf:DG_Ku}
    \mathcal{K}_\gamma(\eta) &=& a_i + \sum_{k=1}^m \left\{ 2 b_{ik} \mu_k + \frac{[2\bP^\top \bH^\top \bB]_{ki}d_k \eta}{1 - 2 \lambda_k \eta}\right\}
 \end{eqnarray*}
 which can be shown to be analytic at $\hat\eta$.
 Consequently, we have
 \begin{equation*}\label{cgf:DG_dKu}
    \frac{\partial \mathcal{K}_\gamma(\eta)}{\partial \eta} = \sum_{k=1}^m \frac{[2\bP^\top \bH^\top \bB]_{ki}d_k}{(1 - 2 \lambda_k \eta)^2} \mbox{ and } \frac{\partial^2 \mathcal{K}_\gamma(\eta)}{\partial \eta^2} = \sum_{k=1}^m \frac{4 \lambda_k [2\bP^\top \bH^\top \bB]_{ki}d_k}{(1 - 2 \lambda_k \eta)^3}.
 \end{equation*}

 Now, we are ready to compute the VaR sensitivity and CVaR sensitivity with respect to $\mu_i$ as
 $$
 \frac{\partial v_\alpha(Y)}{\partial \mu_i} = \mathsf{E}\left[\left.\frac{\partial Y}{\partial \mu_i}  \right| Y = v_\alpha(Y)\right] \mbox{ and } \frac{\partial c_\alpha(Y)}{\partial \mu_i} = \mathsf{E}\left[\left.\frac{\partial Y}{\partial \mu_i}  \right| Y \geq v_\alpha(Y)\right].
 $$
 All the assumptions in \cite{Hong:08} and \cite{Hong Liu:09} are satisfied in this setting.
 Any root-finding algorithm can be applied to locate the unique saddlepoint $\hat\eta$. Once we find $\hat{\eta}$ with the CGF (\ref{cgf:DG_Y}) of $Y$, we are able to derive saddlepoint approximations of risk sensitivities utilizing Theorems \ref{thm1} and \ref{thm2}, as summarized in the following theorem.

 \begin{thm}\label{prop:DG}
    The VaR and CVaR sensitivities with respect to $\mu_i$, the mean of a risk factor, of a delta-gamma portfolio loss $Y$ in (\ref{DG_Y}) are approximated via saddlepoint techniques by
    \begin{eqnarray*}
        \frac{\partial v_\alpha(Y)}{\partial \mu_i} &=& a_i + \sum_{k=1}^m \left\{ 2 b_{ik} \mu_k + \frac{[2\bP^\top \bH^\top \bB]_{ki}d_k}{\left(1 - 2 \lambda_k \hat{\eta}\right)^3} \Bigg[ \hat{\eta} (1 - 2 \lambda_k \hat{\eta})^2 \right. \\
        && \left. \left. + \frac{\cK^{(3)}_Y(\hat{\eta})(1 - 2 \lambda_k \hat{\eta}) -  4 \lambda_k \cK''_Y(\hat{\eta})}{ 2n \cK''_Y(\hat{\eta})^2 + \frac{1}{4} \cK^{(4)}_Y(\hat{\eta}) - \frac{5}{12} \cK^{(3)}_Y(\hat{\eta})^2 /\cK''_Y(\hat{\eta})} \right] \right\}, \nonumber
    \end{eqnarray*}
    and
    \begin{eqnarray*}
        \frac{\partial c_\alpha(Y)}{\partial \mu_i} &=& a_i + \sum_{k=1}^n \left\{ 2 b_{ik} \mu_k + \frac{\phi(\sqrt{n} \hat{\omega})}{\sqrt{n}\hat{z}(1-\alpha)} \cdot \frac{[2\bP^\top \bH^\top \bB]_{ki}d_k}{1 - 2 \lambda_k \hat{\eta}} \right. \times \nonumber \\
        && \left[ \hat{\eta} + \frac{1}{n} \left(  \left(\frac{1}{8}\hat{\rho}_4 - \frac{5}{24}\hat{\rho}_3^2 - \frac{\hat{\rho}_3}{2\hat{z}} - \frac{1}{\hat{z}^2} \right)\hat{\eta}  \right.\right. + \\
        && \left.\left.\left. \frac{1}{\sqrt{\mathcal{K}_{Y}''(\hat{\eta})}(1 - 2 \lambda_k \hat{\eta})} \left( \frac{\hat{\rho}_3}{2} + \frac{1}{\hat{z}} \right) - \frac{4 \lambda_k }{2 \mathcal{K}_{Y}''(\hat{\eta})(1 - 2 \lambda_k \hat{\eta})^2} \right) \right] \right\},
    \end{eqnarray*}
    respectively. The saddlepoint $\hat{\eta}$ is the unique solution of
    \begin{equation*}
        \sum_{i=1}^n \frac{\lambda_i (1 - 2 \lambda_i \eta) + d_i^2 (1 - \lambda_i \eta)\eta}{(1 - 2 \lambda_i \eta)^2} = v_\alpha(Y) - c.
    \end{equation*}
    Here, $\hat{\omega} = \sqrt{2\left(\hat{\eta}a - \mathcal{K}_{Y}(\hat{\eta})\right)}$, $\hat{z} = \hat{\eta}\sqrt{\mathcal{K}_{Y}''(\hat{\eta})}$, and the standardized cumulants are $\hat{\rho}_3 = \cK^{(3)}_Y(\hat{\eta}) / \cK''_Y(\hat{\eta})^{3/2}$, $\hat{\rho}_4 = \cK^{(4)}_Y(\hat{\eta}) / \cK''_Y(\hat{\eta})^{2}$.
 \end{thm}

 To check numerical performances of our expansions, let us take the same example as appeared in Section 5.1 in \cite{Hong Liu:09}.
 Let $f(\mathbf{0}) = 0.3 $, $\ba = [0.8, 1.5]^\top$ and $\bB = \left[\begin{array}{cc}
                                  1.2 & 0.6 \\
                                  0.6 & 1.5
                                \end{array}\right].$ The risk factor $\bX$ follows $\mathcal{N}(\bmu, \bSigma)$ with $\bmu = [0.01, 0.03]^\top$ and $\bSigma = \left[\begin{array}{cc}
                                  0.02& 0.01 \\
                                  0.01 & 0.02
                                \end{array}\right].$
 For comparison, we compute IPA estimates using $10^7$ observations of $Y$ with the batch size $2000$.
 An asymptotically valid $100(1-\beta)\%$ confidence interval of the VaR sensitivity is also reported, see Section 6 in \cite{Hong:08}.

 Figure \ref{fig:DG-VaR} (i) depicts the VaR sensitivities with respect to $\mu_1$ varying $\alpha$ from $0.9$ to $0.99$.
 As in Section \ref{sec:ch5-app1}, two saddlepoint approximations are given based on VaR estimates using simulation and saddlepoint techniques; We denote them by `SPA from MC-VaR' and `SPA from SPA-VaR', respectively.
 \begin{figure}[h]
\centering
\renewcommand{\thesubfigure}{(i)}
\subfigure[][]{\includegraphics[width=0.49\linewidth, height=.36\linewidth]{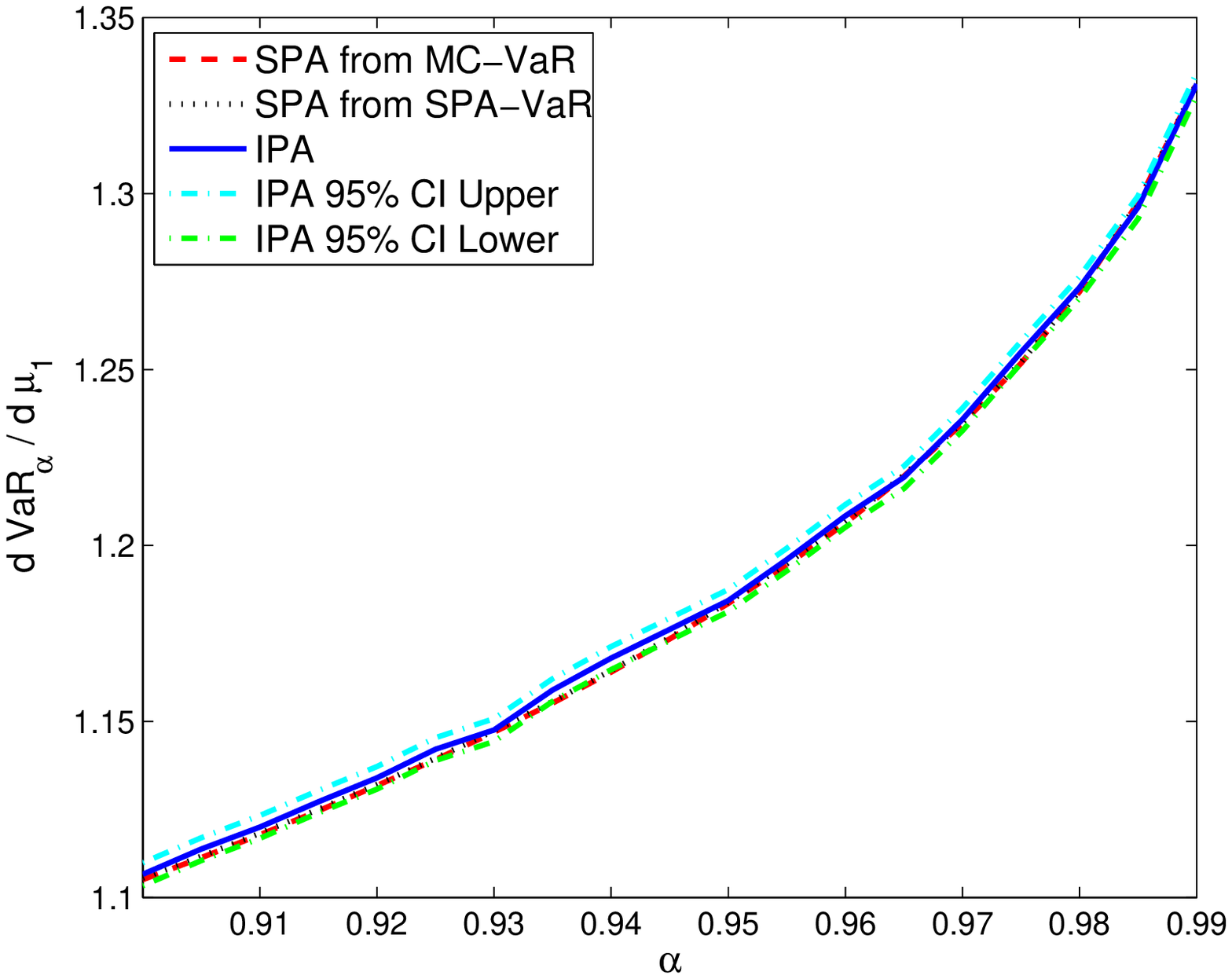}}
\renewcommand{\thesubfigure}{(ii)}
\subfigure[][]{\includegraphics[width=0.49\linewidth, height=.36\linewidth]{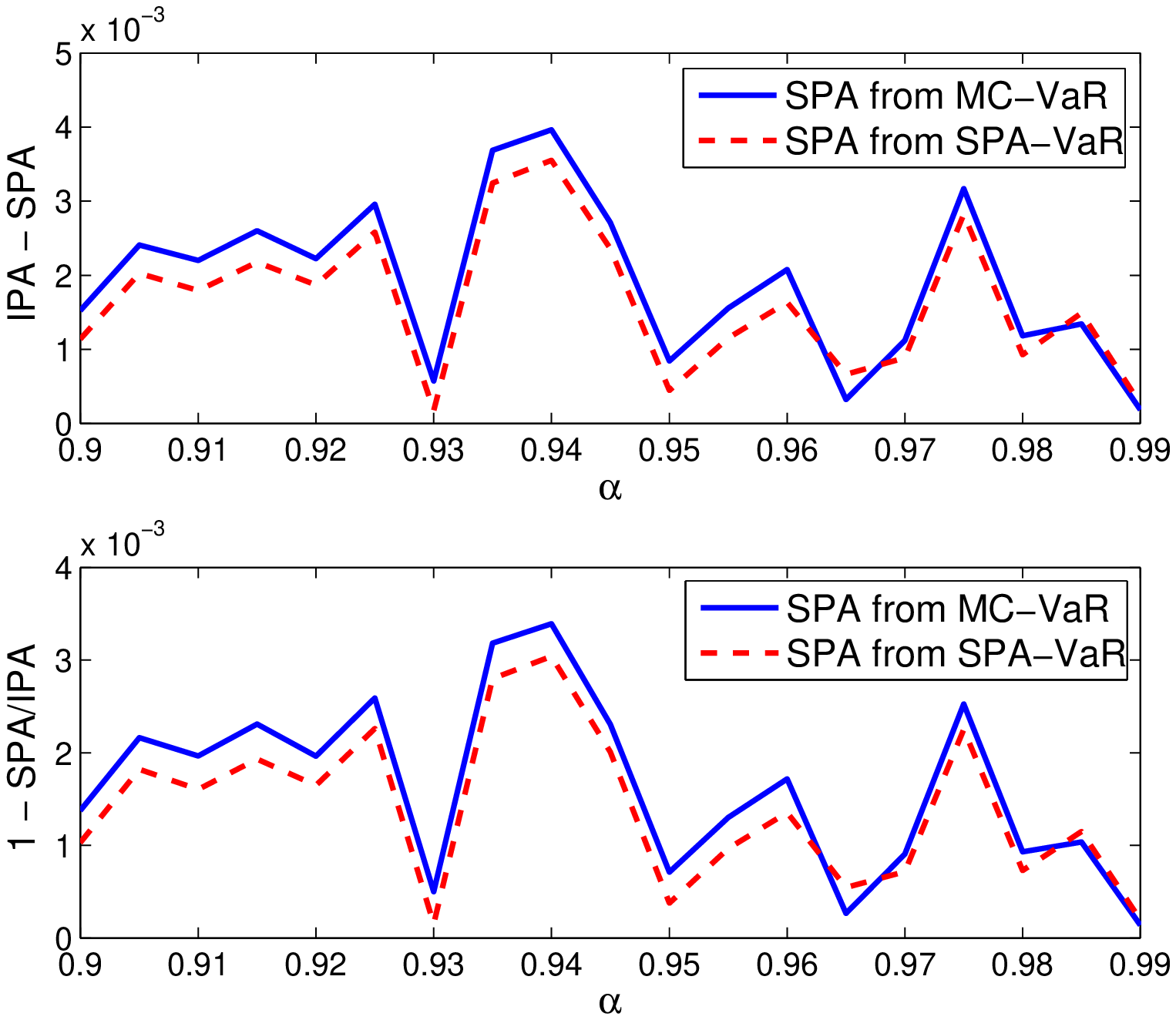}}
\caption{(i) VaR sensitivity with respect to $\mu_1$ over $\alpha$ using SPAs (red, black) and IPA estimator (blue) and (ii) the estimated differences and relative differences of our SPAs to IPA estimator.}\label{fig:DG-VaR}
\centering
\end{figure}
 The blue curve is the IPA estimates together with 95\% confidence interval;  `CI Upper' for the upper bound and `CI Lower' for the lower bound of the interval. The batch size $k$ has been chosen to make the sample variance reasonably small, specifically, $0.0055$. The differences and relative differences of SPA based estimates compared to IPA estimates are shown in Figure \ref{fig:DG-VaR}(ii).
 This figure tells us that Theorem \ref{prop:DG} provides a highly accurate approximation to the sensitivity of VaR regardless of whether we use saddlepoint methods or Monte Carlo simulation for the estimation of $v_\alpha(Y)$.
 For example, the averaged relative difference for SPA from MC-VaR is reported as $1.6462 \times 10^{-3}$. The average (relative) difference between the two VaR sensitivities from MC-VaR and SPA-VaR is even smaller as $3.4591 \times 10^{-4}$ ($2.9550 \times 10^{-4}$).

 Figure \ref{fig:DG-ES} (i) plots the CVaR sensitivities with respect to $\mu_1$ varying $\alpha$ from $0.9$ to $0.99$.
 Similarly as above, we estimate $v_\alpha(Y)$ by IPA or saddlepoint methods, denoting the results by (red) SPA from MC-VaR, (black) SPA from SPA-VaR.
 We also draw IPA estimates as well as interval estimates. Part (ii) of the figure shows the errors and the relative differences of saddlepoint approximations compared to IPA estimates.
 As seen from Figure \ref{fig:DG-ES}, we again see that the expansion in Theorem \ref{prop:DG} gives very fast and accurate results. We, however, note that there are larger differences between the two SPA based estimates (MC-VaR vs. SPA-VaR) than in the case for the VaR sensitivity. The average difference between SPA from MC-VaR and the IPA estimates is $7.2 \times 10^{-4}$ whereas SPA from SPA-VaR gives $1.58 \times 10^{-3}$.

\begin{figure}[h]
\centering
\renewcommand{\thesubfigure}{(i)}
\subfigure[][]{\includegraphics[width=0.49\linewidth, height=.36\linewidth]{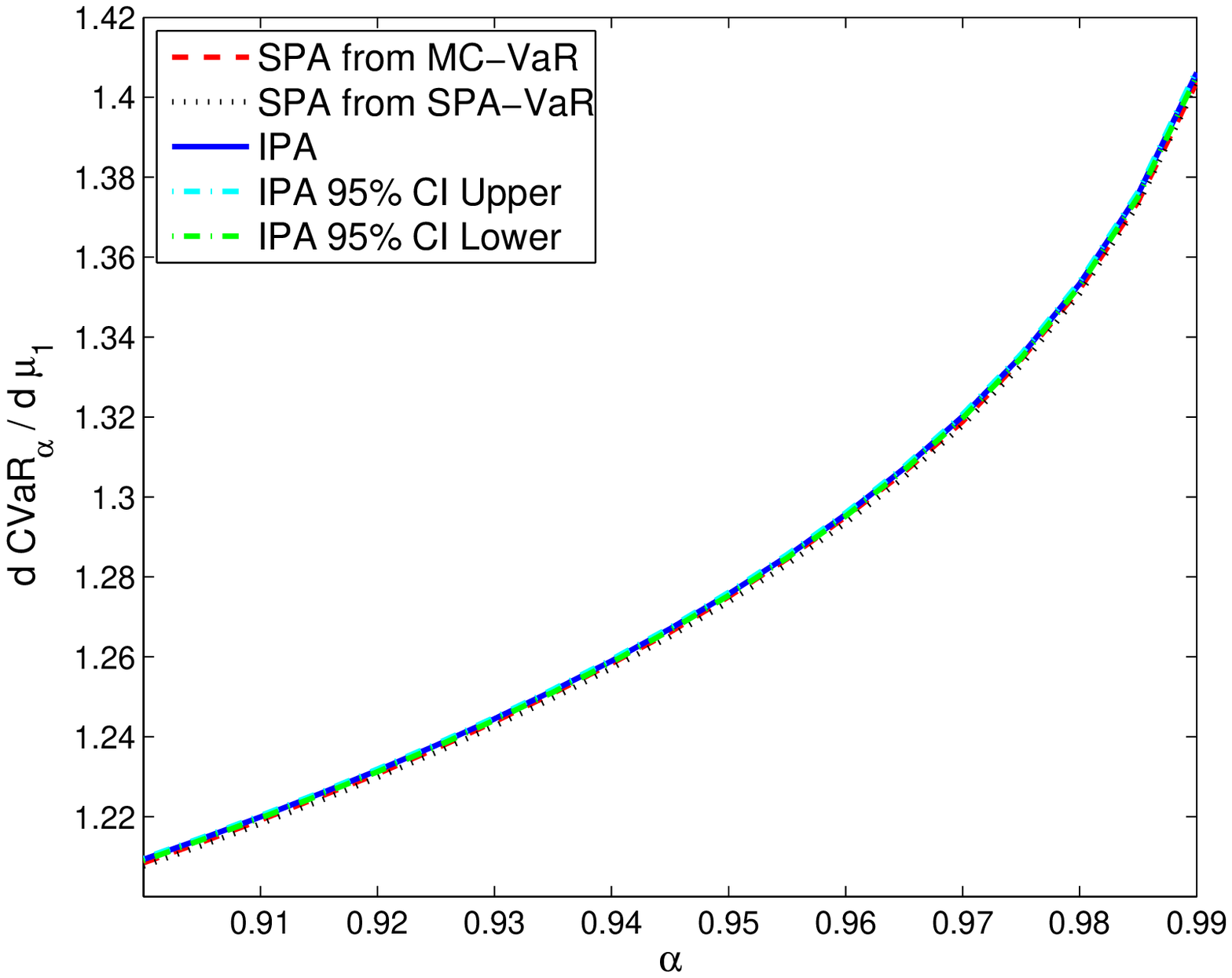}}
\renewcommand{\thesubfigure}{(ii)}
\subfigure[][]{\includegraphics[width=0.49\linewidth, height=.36\linewidth]{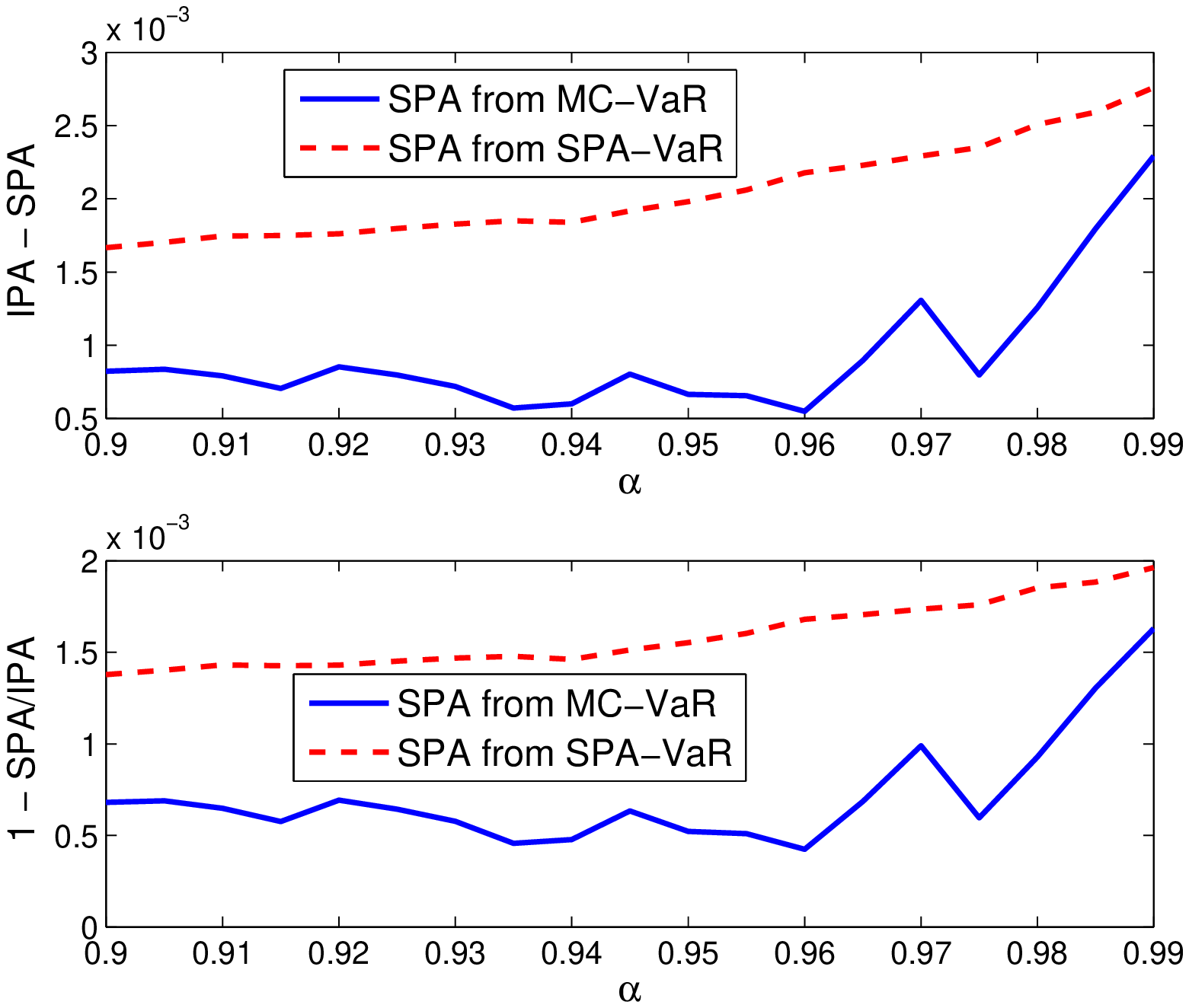}}
\caption{(i) CVaR sensitivity with respect to $\mu_1$ over $\alpha$ using SPAs (red, black) and IPA estimator (blue) and (ii) the estimated differences and relative differences of our SPAs to IPA estimator.}\label{fig:DG-ES}
\centering
\end{figure}

 \subsection{Option sensitivity}\label{sec:ch5-app3}

 Computing sensitivities or {\it greeks} of an option price with respect to market parameters is another important application in financial risk management.
 An option price is typically expressed in terms of the expectation of a payoff functional of underlying asset prices under the risk neutral measure. And its sensitivities can also be expressed as expectations of derivatives of the payoff functional.  For instance, Theorem 1 in \cite{Hong Liu:11} proves that under certain technical conditions the sensitivity of $p(\theta) = \mathsf{E}[g(S) {\bf 1}_{[h(S) \geq 0]}]$ with respect to a parameter $\theta$ is given by
 \begin{equation*}
    \frac{\partial p(\theta)}{\partial \theta} = \mathsf{E}\left[ \frac{\partial g(S)}{\partial \theta} {\bf 1}_{[h(S) \geq 0]} \right] - \left. \frac{\partial }{\partial y} \mathsf{E}\left[ g(S) \frac{\partial h(S)}{\partial \theta} {\bf 1}_{[h(S) \geq y]} \right]\right|_{y=0}
 \end{equation*}
 where $S = \{S(t)\}_{0 \leq t \leq T}$ denotes the underlying asset process.
 This problem has been extensively studied in the literature both by academics and practitioners. Popular methods include finite difference scheme, the pathwise method (equivalent to IPA), the likelihood ratio method, Malliavin calculus, etc. Our objective is to tackle the problem by employing our saddlepoint expansions.

 We choose to work on financial options with two underlying assets and study their sensitivities with respect to volatilities, so called {\it vega}. This is for an illustrative purpose and we note that there are many other possibilities. Furthermore, a bivariate geometric Brownian motion process and an exponential variance gamma model are adopted for the underlying asset processes.

 \subsubsection{Two-asset correlation call option under geometric Brownian motions}

 Suppose that an underlying asset $(S_1(t), S_2(t))$ of an option is a bivariate geometric Brownian motion such that each price process is given by
 $$
 S_i(t) = S_i(0) \exp\left(\left(r_i - \frac{1}{2}\sigma_i^2\right)t + \sigma_i W_i(t)\right)
 $$
 where $W_i$ is a standard Brownian motion with $\mathsf{E}[W_1(t)W_2(t)] = \rho t$ for $i=1,2$ under the risk neutral measure $\mathsf{P}$.
 We consider an option based on $(S_1(t), S_2(t))$ whose price is
 $$
  C = e^{-rT} \mathsf{E}\Big[ (S_1(T) - K)^+ {\bf 1}_{[S_2(T) > H]} \Big].
  $$
 Then the sensitivity of $C$ with respect to $\sigma_1$ can be computed by
 \begin{eqnarray*}
    \frac{\partial C}{\partial \sigma_1} &=& e^{-rT} \mathsf{E}\left[ \frac{\partial S_1(T)}{\partial \sigma_1} {\bf 1}_{[ S_1(T) > K ]}{\bf 1}_{[ S_2(T) > H ]} \right]  \\
    &=& S_1(0)e^{(r_1 - r - \frac{1}{2}\sigma_1^2)T} \left\{ \mathsf{E}\left[ W_1(T)e^{\sigma_1 W_1(T)} {\bf 1}_{[ W_1(T) > k ]}{\bf 1}_{[ W_2(T) > h ]} \right] \right. \\
    && \left.- \sigma_1 T \ \mathsf{E}\left[e^{\sigma_1 W_1(T)} {\bf 1}_{[ W_1(T) > k ]}{\bf 1}_{[ W_2(T) > h ]} \right] \right\}
 \end{eqnarray*}
 where $k = (\log(K/S_1(0)) - (r_1 - \sigma_1^2/2)T)/\sigma_1$ and $h = (\log(H/S_2(0)) - (r_2 - \sigma_2^2/2)T)/\sigma_2$.

 Let $X = W_1(T)$ and $\bY = (Y_1, Y_2) = (W_1(T), W_2(T))$.
 Under $\mathsf{P}$, the CGF of $\bY$ is given by $\mathcal{K}(\eta_1, \eta_2) = T(\eta_1^2/2 + \rho \eta_1\eta_2 + \eta_2^2/2)$. Let $\mathsf{Q}$ be defined by the Radon-Nikodym derivative
 $$
 \frac{\rd \mathsf{Q}}{\rd \mathsf{P}} = \frac{e^{\sigma_1 X}}{\mathsf{E}[e^{\sigma_1 X}]}.
 $$
 It then follows that
  \begin{equation}\label{dC}
    \frac{\partial C}{\partial \sigma_1} = S_1(0)e^{(r_1 - r - \frac{1}{2}\sigma_1^2)T + \cK(\sigma_1,0)} \Bigg\{ \mathsf{E}^{\mathsf{Q}}\Big[ X {\bf 1}_{[ Y_1 > k ]}{\bf 1}_{[ Y_2 > h ]} \Big] - \sigma_1 T \ \mathsf{P}^{\mathsf{Q}}\Big[ Y_1 > k, Y_2 > h  \Big] \Bigg\}.
 \end{equation}
 Thus, we can approximate the expectation under $\mathsf{Q}$ in (\ref{dC}), ${\sf E}^{\sf Q}$,  by Theorem \ref{thm3}. The second term is also approximated by the existing multivariate tail probability approximation and thus we skip its discussion.

 The CGFs of $\bY$ and $(X, \bY)$ under $\mathsf{Q}$ are computed as follows:
 \begin{eqnarray*}
    \mathcal{K}_{\bY}(\eta_1, \eta_2) &=& \mathcal{K}(\sigma_1 + \eta_1, \eta_2) - \mathcal{K}(\sigma_1,0) \;\; \mbox{ and } \\
    \mathcal{K}_{X,\bY}(\gamma, \eta_1, \eta_2) &=& \mathcal{K}(\sigma_1 + \gamma + \eta_1, \eta_2) - \mathcal{K}(\sigma_1,0).
 \end{eqnarray*}
 The saddlepoint of $\mathcal{K}_{\bY}(\eta_1, \eta_2)$ is obtained as
 $$
 (\hat{\eta}_1,  \hat{\eta}_2) = \left(\frac{k - \rho h}{T(1-\rho^2)}, \frac{h - \rho k}{T(1-\rho^2)}\right).
 $$
 Similarly, $\tilde{\eta}_2(\eta_1) = h/T - \rho(\eta_1 + \sigma_1)$ and  $\tilde{\eta}_2(0) = h/T - \rho \sigma_1$.
 The assumption of Theorem \ref{thm3} is satisfied since $\mathcal{K}_{\gamma}(\eta_1, \eta_2) = T(\eta_1 + \sigma_1 + \rho \eta_2)$ is analytic at $(\hat{\eta}_1,\hat{\eta}_2)$.
 All the variables that appear in Theorem \ref{thm3} can be explicitly computed in this setting.

 As the saddlepoint equation is solved analytically and the CGFs under consideration are at most quadratic functions, the relations among the variables $\boldeta$, $\bomega$, and $\bv$ are tractable. Therefore, we can easily compute $\tilde{\omega}_2(\omega_1) = \rho \omega_1/ \sqrt{1-\rho^2}$ so that $F = 0$. In addition,
 $\tilde{\mathcal{K}}_\gamma(v_1,v_2) = T\sigma_1 + \sqrt{T}v_1/\sqrt{1-\rho^2} + \rho \sqrt{T}v_2$ by employing the inverse functions of $v_1(\eta_1)$ and $v_2(\eta_1, \eta_2)$, which are $k_1(v_1) = \sqrt{T/(1-\rho^2)}$ and $k_2(v_2) = \rho \sqrt{T}$. With $n=1$, $\hat{x} = k/\sqrt{T} - \sqrt{T} \sigma_1$ and $\hat{y} = h/\sqrt{T} - \sqrt{T} \rho \sigma_1$.
 Finally, we arrive at the following saddlepoint expansion:
 \begin{equation}\label{EQ_GBM}
    \mathsf{E}^{\mathsf{Q}}\Big[ X {\bf 1}_{[ Y_1 > k ]}{\bf 1}_{[ Y_2 > h ]} \Big] \approx \sigma_1 T \bar{\Phi}(\hat{x}, \hat{y}, \rho) + \sqrt{T} \left[ \phi(\hat{x})\bar{\Phi}\left(\frac{\hat{y} - \rho \hat{x}}{\sqrt{1-\rho^2}}\right) + \rho\phi(\hat{y})\bar{\Phi}\left(\frac{\hat{x} - \rho \hat{y}}{\sqrt{1-\rho^2}}\right) \right].
 \end{equation}
 The true value of $\mathsf{E}^{\mathsf{Q}}\Big[ X {\bf 1}_{[ Y_1 > k ]}{\bf 1}_{[ Y_2 > h ]} \Big]$ can be computed as $\bY$ follows a bivariate normal distribution $\mathcal{N}(\sigma_1 T , \rho \sigma_1 T, T, T, \rho)$ under $\mathsf{Q}$, namely,
 \begin{equation}\label{EQ_GBM_true}
    \mathsf{E}^{\mathsf{Q}}\left[ X {\bf 1}_{[ Y_1 > k ]}{\bf 1}_{[ Y_2 > h ]} \right] = \sigma_1 T \bar{\Phi}(\hat{x}, \hat{y}, \rho) + \sqrt{T} \int_{\hat{x}}^\infty \int_{\hat{y}}^\infty y_1 \phi_{\rho}(y_1, y_2) dy_1 dy_2
 \end{equation}
 where $\phi_{\rho}(y_1, y_2)$ is a joint PDF of $\mathcal{N}(0 ,0, 1, 1, \rho)$.
 And it turns out that (\ref{EQ_GBM}) and (\ref{EQ_GBM_true}) coincide.

 \subsubsection{Exchange option under exponential variance gamma models}

 In the second example, we consider an exchange option whose risk neutral valuation formula is given by
 $$
 C = e^{-rT} \mathsf{E}\Big[\big(S_1(T) - S_2(T)\big)^+\Big]
 $$
 based on two assets $(S_1(t), S_2(t))$. Each $S_i(t)$ is assumed to be an exponential variance gamma process, e.g., $S_i(t) = S_i(0)\exp\big(r_i t + \sigma_iX_i(t)\big)$ where
 $X_i(t)$ is an independent variance gamma process.
 The CGF of $X_i(T)$ under the risk neutral measure $\mathsf{P}$ is
 $$
 \mathcal{K}_i(\gamma) =  -\frac{T}{v_i} \log \left( 1 - \theta_i v_i \gamma - \frac12\kappa_i v_i\gamma^2 \right)
 $$
 for the parameter set $(\theta_i, \kappa_i , v_i)$. Note that $X_i(t)$ can be interpreted as a time-changed Brownian motion such that $X_i(t) = \theta_i G_i(t) + \kappa_i W_i(G_i(t))$ where $G_i(t)$ is a gamma process independent of $W_i$ with unit drift and volatility $v_i$.
 We also denote the CGF of $(X_1(T), X_2(T))$ under ${\sf P}$ by $\cK(\eta_1, \eta_2)$.

 We are interested in the sensitivity of the option price $C$ with respect to $\sigma_1$.
 It can be computed by
 \begin{eqnarray}\label{dC2}
    \frac{\partial C}{\partial \sigma_1} &=& S_1(0) \mathsf{E}\left[ \frac{\partial S_1(T)}{\partial \sigma_1} {\bf 1}_{[ S_1(T) > S_2(T) ]} \right] \nonumber \\
    &=& S_1(0)e^{(r_1 - r)T + \mathcal{K}_1(\sigma_1)} \mathsf{E}^{\mathsf{Q}} \Big[ X_1(T) {\bf 1}_{[ \sigma_1 X_1(T) - \sigma_2 X_2(T) > k ]}\Big]
 \end{eqnarray}
 where $k = \log(S_2(0)/S_1(0)) +(r_2 - r_1)T$
 and $\mathsf{Q}$ is again defined by the Radon-Nikodym derivative
 $ {\rd \mathsf{Q}}/{\rd \mathsf{P}} = {e^{\sigma_1 X_1(T)}}/{\mathsf{E}[e^{\sigma_1 X_1(T)}]}.$
 Take $X = X_1(T)$ and $Y = \sigma_1 X_1(T) - \sigma_2 X_2(T)$.
 The CGF of $X$ is $\mathcal{K}_X(\gamma) = \mathcal{K}_1(\gamma + \sigma_1) - \mathcal{K}_1(\sigma_1)$; the CGF of $Y$ is $\mathcal{K}_Y(\eta) = \mathcal{K}((1+\eta)\sigma_1, -\eta\sigma_2) - \mathcal{K}(\sigma_1,0)$; and the joint CGF of $(X,Y)$ then is obtained by
 $\mathcal{K}_{X,Y}(\gamma,\eta) = \mathcal{K}(\sigma_1 + \gamma + \sigma_1 \eta, -\sigma_2 \eta) - \mathcal{K}(\sigma_1,0)$ under $\mathsf{Q}$.

 The convergence domain of the above CGFs contain zero. And the saddlepoint $\hat{\eta}$ of $Y$ is the solution of a polynomial equation of degree four, which can be numerically found by the Newton-Raphson method.
 Moreover, $\mathcal{K}_{\gamma}(\eta)$ is an analytic function in the convergence domain of $\mathcal{K}_Y$ and is given by
 $$ \mathcal{K}_{\gamma}(\eta) = \frac{T(\theta_1 + \kappa_1 v_1 \sigma_1(1+\eta))}{1 - \theta_1 v_1 \sigma_1(1+\eta) - \kappa_1 v_1\sigma_1^2(1+\eta)^2/2}.$$
 Then by applying the saddlepoint formula (\ref{rem:spa}) in Remark \ref{rem:integralspa}, we can finally compute $$\mathsf{E}^{\mathsf{Q}} \left[ X_1(T) {\bf 1}_{[ \sigma_1 X_1(T) - \sigma_2 X_2(T) > k ]}\right]$$ in (\ref{dC2}).
 We omit the details of this computation due to its complexity.

 Figure \ref{fig:ET} shows numerical results of the sensitivity of $C$ with respect to $\sigma_1$ under the parameter set given in Table \ref{tb:ET} with $\theta_i = 0$ for brevity.
 IPA estimates are obtained based on  $10^6$ simulated samples of the variance gamma processes under $\mathsf{P}$.
 The average of the estimated relative differences of two approaches is reported as $1.5\times 10^{-3}$, which also shows great performance of the developed approximations.

 \begin{table}[h]
 \centering
 \begin{tabular}{ccccccccccc}
 \hline
 $S_1(0)$ & $S_2(0)$ & $T$ & $r$ & $r_1$ & $r_2$ & $\sigma_2$ & $v_1$ & $v_2$ & $\kappa_1$ & $\kappa_2$ \\ \hline
 90       & 100      & 1   & 0.02 & 0.2   & 0.4   & 1          & 0.2   & 0.25  & 0.1        & 0.32       \\ \hline
 \end{tabular}
 \vs
 \caption{Parameter set of Figure \ref{fig:ET}.}\label{tb:ET}
 \end{table}

 \begin{figure}[h]
 \centering
 \renewcommand{\thesubfigure}{(i)}
 \subfigure[][]{\includegraphics[width=0.49\linewidth, height=.36\linewidth]{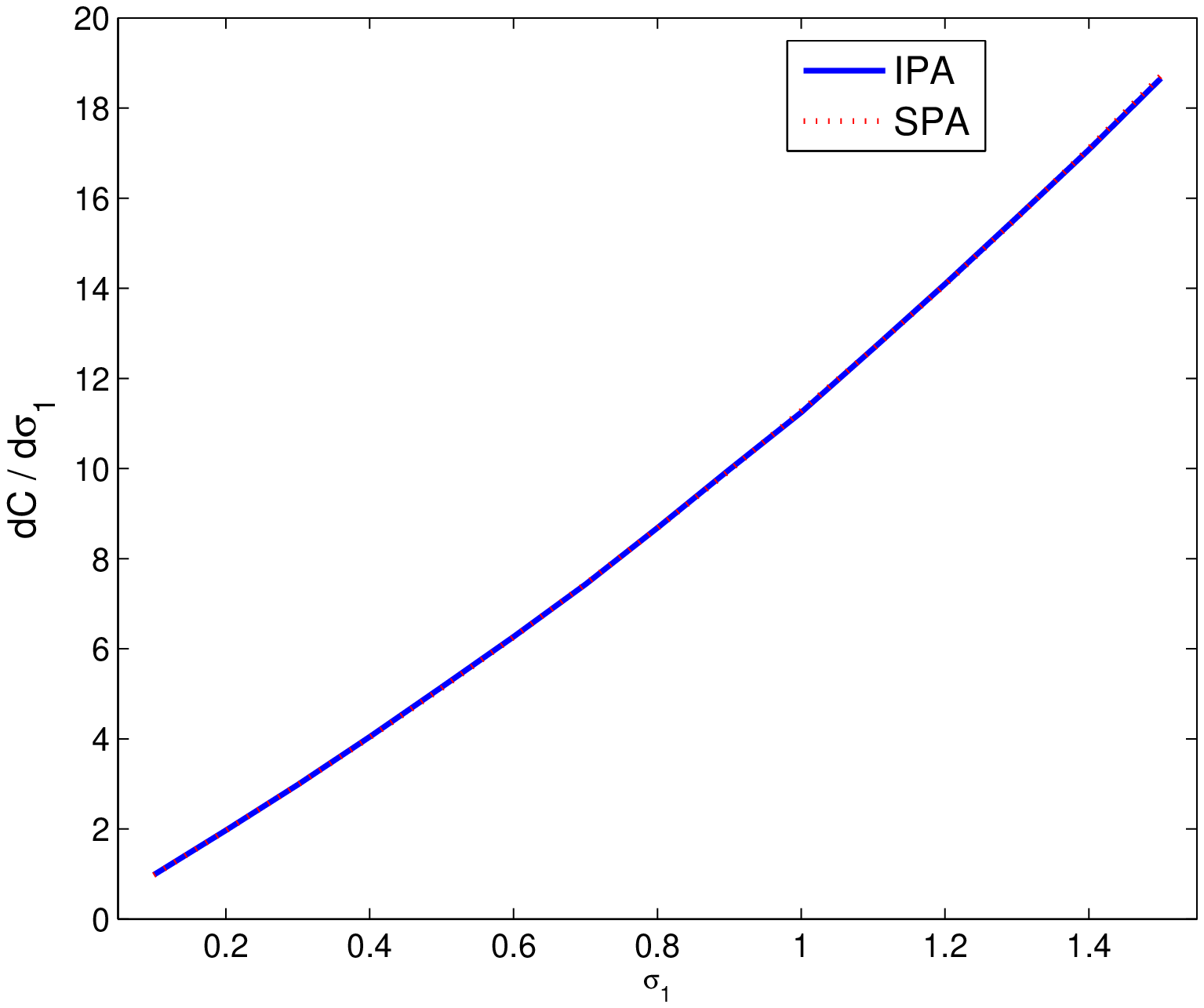}}
 \renewcommand{\thesubfigure}{(ii)}
 \subfigure[][]{\includegraphics[width=0.49\linewidth, height=.36\linewidth]{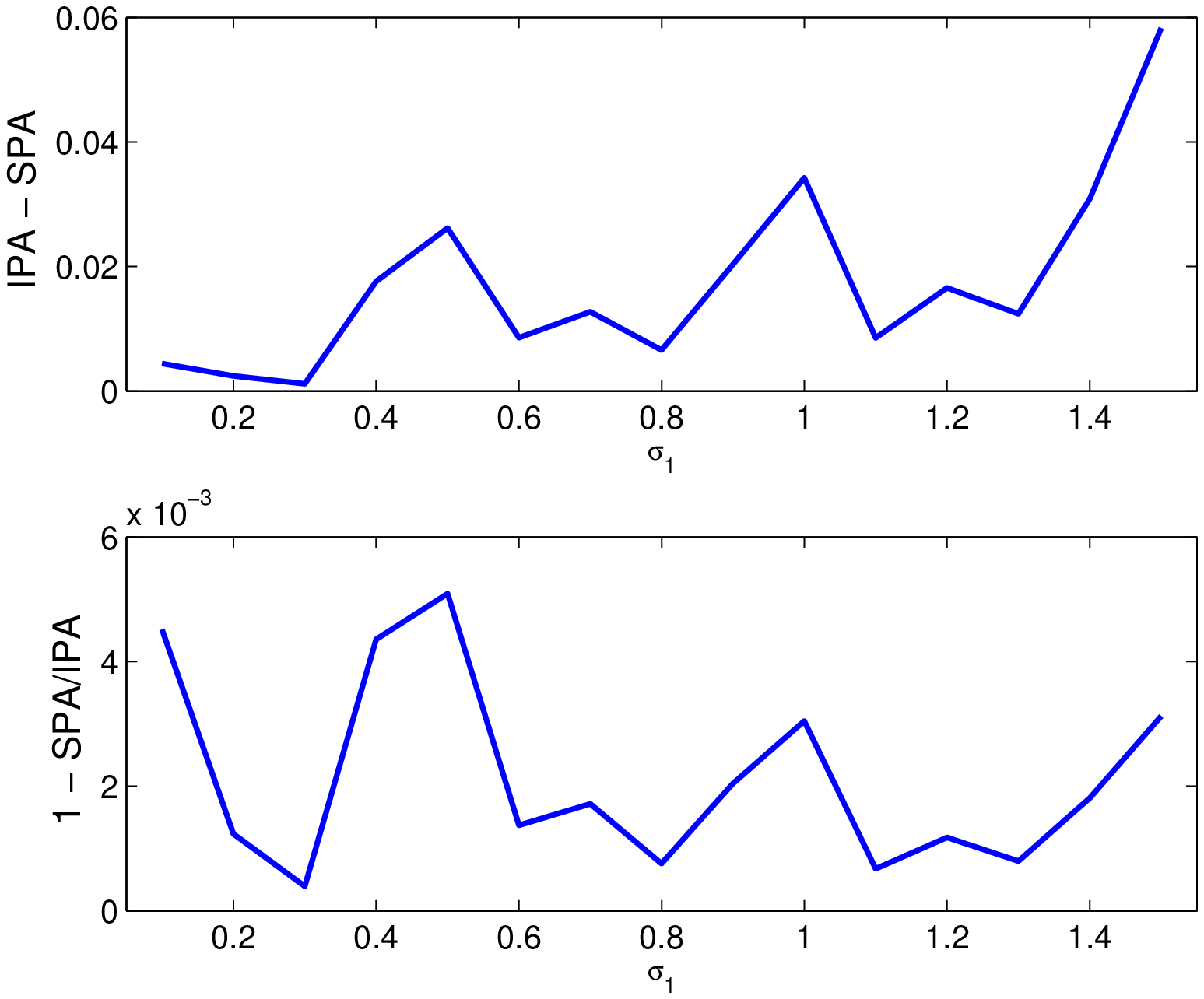}}
 \caption{(i) The sensitivity of an exchange option with respect to $\sigma_1$ via the sadddlepoint method (red) and the IPA estimator (blue) and (ii) the differences and the relative differences of the two values.}\label{fig:ET}
 \centering
 \end{figure}
 
 \section{Conclusion}

 Saddlepoint approximations for $\mathsf{E}[\overline{X} | \overline{Y} = a]$ and $\mathsf{E}[\overline{X} | \overline{Y} \geq a]$ were derived for the sample mean of a continuous bivariate random vector $(X,Y)$ whose joint moment generating function is known. The extensions of the approximations to the case of a random vector $\bY$ were also investigated.
 The newly developed expansions were applied to several problems associated with risk measures and financial options.
 We specifically focused on risk contributions of asset portfolios and risk sensitivities of delta-gamma portfolios.
 Sensitivities of an option based on two assets with respect to a market parameter were also computed via the proposed saddlepoint approximations.
 We have performed numerical experiments, showing that the new approximations are not only computationally efficient but also very accurate compared to simulation based estimates.
 As a whole, our developments have broadened the applicability of saddlepoint techniques by providing explicit and accurate approximations to certain conditional expectations.

\begin{appendix}

 \section{Proof of Lemma \ref{lem:inversion}}\label{app:proof of lemma1}

 We first prove the second inversion formula (\ref{lem:cdf-inversion}).
 Suppose that $X$ has a non-negative lower bound.
    $$\mathsf{E}[X {\bf 1}_{[\bY \geq \ba]}]  = \mathsf{E}[X] \int_{[\bY \geq \ba]} \int_0^\infty \frac{x}{\mathsf{E}[X]} f_{X, \bY}(x,\by) \rd x \rd \by = \mathsf{E}[X] \cdot \mathsf{P}_h[\bY \geq \ba],$$
 where the density of $\bY$ under $\mathsf{P}_h$ is $h(\by) = \int_{0}^{\infty} (x/\mathsf{E}[X]) f_{X, \bY}(x,\by) \rd x$.

 The MGF of $\bY$ under $\mathsf{P}_h$ is then
    \begin{eqnarray*}
        \mathcal{M}_h(\boldeta) &=& \int_{\mathbb{R}^d} e^{\by^\top \boldeta} h(\by) \rd\by\\
        &=& \frac{\mathcal{M}_\mathbf{Y}(\boldeta)}{\mathsf{E}[X]} \int_{\mathbb{R}^d}  \int_{0}^{\infty} x \frac{e^{\by^\top \boldeta}}{\mathcal{M}_\bY(\boldeta)} f_{X, \bY}(x,\by)  \rd x \rd \by \\
        &=& \frac{\mathcal{M}_\mathbf{Y}(\boldeta)}{\mathsf{E}[X]} \int_{0}^{\infty} x g(x)\rd x = \frac{\mathcal{M}_\bY(\boldeta)}{\mathsf{E}[X]}\mathsf{E}_g[X],
    \end{eqnarray*}
 where $\mathcal{M}_\bY(\boldeta)$ denotes the MGF of $\bY$ under $\mathsf{P}$, $g(x) = \int_{\mathbb{R}^d} (e^{\by^\top \boldeta}/\mathcal{M}_\bY(\boldeta)) f_{X, \mathbf{Y}}(x,\by)\rd \by  $, and $\mathsf{E}_g$ denotes the integration under the new probability $\mathsf{P}_g$ having the density $g(x)$. The third equality holds by the Fubini theorem due to the non-negativity of the integrand.

 On the other hand, the MGF of $X$ under $\mathsf{P}_g$ can also be computed as
    $$
    \mathcal{M}_g(\gamma) = \int_{0}^{\infty} e^{\gamma x} g(x) \rd x = \frac{\mathcal{M}_{X, \bY}(\gamma, \boldeta)}{\mathcal{M}_\bY(\boldeta)}.
    $$
 Therefore,
    $$
    \mathsf{E}_g[X] = \mathcal{M}_g'(0) = \frac{1}{\mathcal{M}_\bY(\boldeta)} \cdot  \left.\frac{\partial}{\partial \gamma} \mathcal{M}_{X, \bY}(\gamma, \boldeta) \right|_{\gamma=0}
    $$
 so that
    $$
    \mathcal{M}_h(\boldeta) = \frac{1}{\mathsf{E}[X]} \cdot \left.\frac{\partial}{\partial \gamma} \mathcal{M}_{X, \bY}(\gamma, \boldeta) \right|_{\gamma=0}.
    $$
 Thus we obtain the CGF $\mathcal{K}_h(\boldeta)$ under $\mathsf{P}_h$ as
    $$
    \mathcal{K}_h(\boldeta) = \mathcal{K}_\bY(\boldeta) + \log\left.\frac{\partial}{\partial \gamma} \mathcal{K}_{X, \bY}(\gamma, \boldeta) \right|_{\gamma=0} - \log \mathsf{E}[X],
    $$
 since $\partial \mathcal{M}_{X, \bY}(\gamma, \boldeta)/\partial \gamma = \partial \mathcal{K}_{X, \bY}(\gamma, \boldeta)/\partial \gamma  \cdot \mathcal{M}_{X, \bY}(\gamma, \boldeta)$.
 By substituting $\mathcal{K}_h(\boldeta)$ to the inversion formula (\ref{multi-cdf-inversion}), i.e.
    $$
    \mathsf{P}_h[\bY\geq \ba] = \left(\frac{1}{2\pi i}\right)^d \int_{\btau-i\infty}^{\btau + i\infty} \frac{\exp \left( \mathcal{K}_{h}(\boldeta) - \ba^\top \boldeta \right)}{\prod_{j=1}^d \eta_j} \rd \boldeta,
    $$
    we have the desired result.

 In the case that $X$ has a negative lower bound $-B$ with $B>0$, define $Z = X + B$ so that $Z$ has a non-negative lower bound. Then the marginal CGF of $Z$ is $\mathcal{K}_Z(\gamma) = \mathcal{K}_X(\gamma) +B\gamma$ and the joint CGF of $Z$ is $\mathcal{K}_{Z, \mathbf{Y}}(\gamma, \boldeta) = \mathcal{K}_{X, \mathbf{Y}}(\gamma, \boldeta) + B\gamma$ where $\mathcal{K}_X(\gamma)$ denotes the CGF of $X$. Note that
    \begin{eqnarray*}
        \mathsf{E}[X {\bf 1}_{[\bY \geq \ba]}] &=& \mathsf{E}[ (Z-B) {\bf 1}_{[\bY \geq \ba]}] = \mathsf{E}[Z {\bf 1}_{[\bY \geq \ba]}] - B \mathsf{P}[\bY \geq \ba] \\
        &=& \left(\frac{1}{2\pi i}\right)^d \int_{\btau - i\infty}^{\btau + i\infty} \left\{  \left. \frac{\partial}{\partial \gamma} \mathcal{K}_{X, \bY}(\gamma, \boldeta) \right|_{\gamma=0} + B \right\} \\
        && \times \frac{\exp \left( \mathcal{K}_{\bY}(\boldeta) - \ba^\top \boldeta \right)}{\prod_{j=1}^d \eta_j} \rd \boldeta - B \mathsf{P}[\bY \geq \ba]
    \end{eqnarray*}
 from the result for the non-negative case.
 This immediately leads to (\ref{lem:cdf-inversion}).

 Finally, for an unbounded $X$, we take $X_C = \max(X,C)$ where $C$ is a constant. The assumption imposed on the MGF of $(X, \bY)$ implies that the MGF ${\mathcal M}_{X_C,\bY}$ also exists in an open neighborhood of the origin. Since $X_C$ is bounded from below,
    \begin{equation}\label{lem:MCT}
        \mathsf{E}[X_C {\bf 1}_{[\bY \geq \ba]}] = \left(\frac{1}{2\pi i}\right)^d \int_{\btau-i\infty}^{\btau+i\infty} \left. \frac{\partial}{\partial \gamma}  \mathcal{M}_{X_C,\bY}(\gamma,\boldeta) \right|_{\gamma=0} \frac{\exp(- \ba^\top \boldeta)}{\boldeta} \rd\boldeta.
    \end{equation}
 But, we have
    \begin{eqnarray*}
        \lefteqn{\frac{\partial}{\partial \gamma}  \mathcal{M}_{X_C,\bY}(\gamma,\boldeta)} \\
        &=& \int_{\mathbb{R}^d}  \int_{-\infty}^{\infty} (x\vee C)e^{\gamma(x\vee C) + \boldeta^\top \by} f_{X,\bY}(x,\by)\rd x \rd\by \\
        &=& \int_{\mathbb{R}^d}   \int_{-\infty}^C C e^{\gamma C + \boldeta^\top \by} f_{X,\bY}(x,\by)\rd x \rd\by  + \int_{\mathbb{R}^d}  \int_{C}^{\infty} x e^{\gamma x + \boldeta^\top \by} f_{X,\bY}(x,\by)\rd x \rd\by \\
        &=& \frac{\partial}{\partial \gamma} \mathcal{M}_{X, \bY}(\gamma, \boldeta) + \int_{\mathbb{R}^d} \int_{-\infty}^C \left\{ C e^{\gamma C} - xe^{\gamma x} \right\} e^{\boldeta^\top \by} f_{X,\bY}(x,\by)\rd x \rd\by.
    \end{eqnarray*}
 The change of integration and differentiation in the first equality is justified by the continuity of the integrand.
 Thus, with $\gamma =0$ it decreases monotonically as $C$ decreases, and as $C\rightarrow -\infty$ we have
    $$ \left. \frac{\partial}{\partial \gamma}  \mathcal{M}_{X_C,\bY}(\gamma,\boldeta)\right|_{\gamma=0} \rightarrow \left. \frac{\partial}{\partial \gamma}  \mathcal{M}_{X,\bY}(\gamma,\boldeta)\right|_{\gamma=0}. $$
 Since $X_C$ converges to $X$ almost surely and monotonically as $C\rightarrow-\infty$, we obtain (\ref{lem:cdf-inversion}) by applying the monotone convergence theorem to both sides of (\ref{lem:MCT}).

 The first formula (\ref{lem:pdf-inversion}) follows similarly with ease. But, we do not need to set $\btau$ to be positive since the inversion formula (\ref{multi-pdf-inversion}) holds for any $\btau$ in a suitable domain, and the convergence for an unbounded case can be proved after a simple adjustment.

 \section{Proof of Theorem \ref{cor:e-thm1}}\label{app:proof of cor}

 We first demonstrates a rescaled and multivariate version of Watson's lemma, Theorem 6.5.2 in \cite{Kolassa:06}.

 \begin{lem}\label{multi-Watsons}
    Suppose that $\theta_j(\bomega)$'s are analytic functions from a domain $\bQ \subset \mathbb{C}^d$ to $\mathbb{C}$ for $0 \leq j \leq k$, and let
    \begin{equation*}\label{def:f_n}
        \vartheta_n(\bomega) = \sum_{j=0}^k \theta_j(\bomega) / n^j.
    \end{equation*}
    Take $\hat{\bomega} \in \bQ$ such that $\hat{\bomega} + i[-\epsilon, \epsilon]^d \subset \bQ$. Then
    \begin{equation*}
        \left(\frac{n}{2\pi}\right)^{d/2} i^{-d} \int_{\hat{\omega}_1 -i\epsilon}^{\hat{\omega}_1 +i\epsilon} \cdots \int_{\hat{\omega}_d -i\epsilon}^{\hat{\omega}_d +i\epsilon} \exp\left[ \frac{n}{2}\sum_{i=1}^d (\omega_i - \hat{\omega}_i)^2\right] \vartheta_n(\bomega) \ d \bomega= \sum_{s=0}^{k-1} A_s n^{-s} + O(n^{-k}),
    \end{equation*}
    where
    \begin{equation*}
        A_s = \sum_{{\mathcal J}_s} \frac{(-2)^{-\sum_{j=1}^d v_j}}{v_1! \cdots v_d!} \left[ \frac{\partial^{2v_1 + \cdots + 2v_d}}{\partial^{2v_1} w_1 \cdots \partial^{2v_d} w_d} \theta_{s-\sum_{j=1}^d v_j}\right]\left( \hat{\bomega}\right)
    \end{equation*}
    for ${\mathcal J}_s = \{ (v_1, \ldots, v_d) \in \mathbb{N}^d \ | v_1, \cdots, v_d \geq 0, \sum_{j=1}^d v_j \leq s \}$.
\end{lem}

 By a change of variable in (\ref{lem:cdf-inversion}) with $(\overline{X}, \overline{\bY})$ and by the closed curve theorem, we write $\mathsf{E}[\overline{X} {\bf 1}_{[\overline{\bY} = \ba]}]$ as
 \begin{eqnarray}
        \mathsf{E}\left[\overline{X} {\bf 1}_{[\overline{\bY} = \ba]}\right] &=& \left( \frac{n}{2\pi} \right)^{d/2} \exp \left( \mathcal{K}_{\bY}(\hat{\boldeta}) -  \hat{\boldeta}^\top \ba \right)  \label{pf:ce1}\\
        &\times&  i^{-d} \left( \frac{n}{2\pi} \right)^{d/2} \int_{\hat{\bomega}-i\infty}^{\hat{\bomega}+i\infty} \frac{n}{2}(\bomega - \hat{\bomega})^\top(\bomega- \hat{\bomega}) \ \mathcal{K}_{\gamma}(\boldeta(\bomega)) \left| \frac{\partial \boldeta}{\partial \bomega} \right| \rd\bomega. \nonumber
    \end{eqnarray}
 We take $\theta_j(\bomega) = 0$ unless $j = 0$ and $\theta_0(\bomega) = \mathcal{K}_{\gamma}(\boldeta(\bomega)) \left| {\partial \boldeta}/{\partial \bomega} \right|$.
 Applying Lemma \ref{multi-Watsons} to (\ref{pf:ce1}) with $k=2$, $A_0 = \theta_0 (\hat{\omega})$ and $A_1 = - \sum_{i=1}^d (\partial^2 \theta_0 / \partial \omega_i^2) (\hat{\bomega})/2$.
 Obtaining $A_1$ only requires to compute the first and second derivatives of $\det[\partial \boldeta / \partial\bomega]$ and $\eta_k$ with respect to $\omega_k$, evaluated at $\hat{\boldeta}$.
 Since computation of coefficients is messy, we here omit the details but report the following formula in \cite{Kolassa:06}:
 \begin{eqnarray*}
      \frac{\partial \boldeta}{\partial \bomega}( \bomega)&=& \frac{\partial \boldeta}{\partial \bomega} (\hat{\bomega}) \left\{ 1 - \frac{1}{3} \sum_m \sum_{i,j,l} \hat{\kappa}^{ijl} \hat{\kappa}_{il}\hat{\eta}_j^m (\omega_m - \hat{\omega}_m) \right.\\
      && + \frac{1}{2} \sum_{m,n} \sum_{o,l} \left[ \frac{1}{4} \sum_{g,h,i,j} \hat{\kappa}^{gij} \hat{\kappa}_{gh} \hat{\kappa}^{hol}\hat{\kappa}_{ij} + \frac{1}{6}\sum_{g,h,i,j} \hat{\kappa}^{gil} \hat{\kappa}_{gh} \hat{\kappa}^{hjo}\hat{\kappa}_{ij} \right.\\
      && \left.\left. - \frac{1}{4}\sum_{i,j} \hat{\kappa}^{ijol} \hat{\kappa}_{ij}\right] \hat{\eta}_o^m \hat{\eta}_l^m (\omega_m - \hat{\omega}_m)(\omega_n - \hat{\omega}_n)\right\} + O(\|\omega - \hat{\omega}\|).
 \end{eqnarray*}
 From this formula, all the desired quantities can be derived in a messy but straightforward manner.

 \section{Proof of Lemma \ref{lem:I2}}\label{app:I2}
 The sum of two integrals $I^{12}+I^{2}$ is expressed as
 \begin{eqnarray}\label{lem:I2-1}
    & & \int_{\hat{\bomega}-i\infty}^{\hat\bomega+i\infty} \frac{\exp \left[ n q(\omega_1,\omega_2) \right]}{(2\pi i)^2} \frac{H^{12}(\eta_1,\eta_2)}{\omega_1 (\omega_2 - \tilde{\omega}_2(\omega_1))} \mathcal{K}_{\gamma}(\eta_1, \eta_2) \rd\bomega \nonumber \\
    &+& \int_{\hat{\bomega}-i\infty}^{\hat{\bomega}+i\infty} \frac{\exp \left[ n q(\omega_1,\omega_2) \right]}{(2\pi i)^2 \ \omega_1} \frac{H^{2}(\eta_1, \eta_2)}{\omega_2 - \tilde{\omega}_2(\omega_1)} \mathcal{K}_{\gamma}(\eta_1, \eta_2) \rd \bomega.
 \end{eqnarray}

 Let $\theta_0^2(\omega_1, \omega_2) = H^{2}(\eta_1, \eta_2) \mathcal{K}_{\gamma}(\eta_1, \eta_2)/\left( \omega_2 - \tilde{\omega}_2(\omega_1)\right)$ as a function of $(\omega_1, \omega_2)$; next we decompose ${\theta_0^2}/{\omega_1}$ as
 $$\frac{\theta_0^2(\omega_1, \omega_2)}{\omega_1} = \frac{\theta_0^2(0, \omega_2)}{\omega_1} + \frac{\theta_0^2(\omega_1, \omega_2) - \theta_0^2(0, \omega_2)}{\omega_1}.$$
 Then (\ref{lem:I2-1}) can be computed as
 \begin{eqnarray*}
    & & \int_{\hat{\bomega}-i\infty}^{\hat{\bomega}+i\infty} \frac{\exp \left[ n q(\omega_1,\omega_2) \right]}{(2\pi i)^2} \left\{ \frac{H^{12}(\eta_1,\eta_2) \mathcal{K}_{\gamma}(\eta_1, \eta_2)}{\omega_1 (\omega_2 - \tilde{\omega}_2(\omega_1))} + \frac{\theta_0^2(\omega_1, \omega_2) - \theta_0^2(0, \omega_2)}{\omega_1} \right\} \rd \bomega  \\
    & & + \frac{1}{2\pi i} \int_{\hat{\omega}_1-i\infty}^{\hat{\omega}_1+i\infty} \exp\left[ n\left( \frac12 \omega_1^2 - \hat{\omega}_1 \omega_1 \right) \right] \frac{1}{\omega_1} \rd \omega_1 \\
    & & \times \frac{1}{2\pi i} \int_{\hat{\omega}_2-i\infty}^{\hat{\omega}_2 +i\infty} \exp\left[ n\left( \frac12 \omega_2^2 - \hat{\omega}_2 \omega_2 \right) \right] \theta_0^2(0, \omega_2) \rd \omega_2 \\
    &=& \frac{1}{2\pi n} \exp[n(\mathcal{K}_{\bY}(\hat{\eta}_1, \hat{\eta}_2) - \hat{\eta}_1 a_1 - \hat{\eta}_2 a_2)] \theta_0^{12}(\hat{\omega}_1, \hat{\omega}_2) \\
    & & + \frac{1}{\sqrt{n}} \bar{\Phi}(\sqrt{n}\hat{\omega}_1) \phi(\sqrt{n}\hat{\omega}_2) \theta_0^2(0,\hat{\omega}_2) + O(n^{-3/2})
 \end{eqnarray*}
 where
 $$
 \theta_0^{12}(\omega_1, \omega_2) = \left\{ \frac{H^{12}(\eta_1,\eta_2) \mathcal{K}_{\gamma}(\eta_1, \eta_2)}{\omega_1 (\omega_2 - \tilde{\omega}_2(\omega_1))} + \frac{\theta_0^2(\omega_1, \omega_2) - \theta_0^2(0, \omega_2)}{\omega_1} \right\}.
 $$
 The last equality holds by applying multivariate and univariate Watson's lemma, namely Lemma \ref{multi-Watsons} and Lemma \ref{Watsons}. This can be done because $\theta_0^{12}(\omega_1, \omega_2)$ and $\theta_0^{2}(\omega_1, \omega_2)$ are analytic near $(\hat{\omega}_1, \hat{\omega}_2)$; then we retain the first order terms only.

 The proof ends by evaluating the meaningful coefficient as
 \begin{eqnarray*}
    \theta_0^2(0,\hat{\omega}_2)
    &=& \mathcal{K}_{\gamma}(0, \eta_2(0,\hat{\omega}_2)) \frac{F(0, \eta_2(0,\hat{\omega}_2)) - F(0,0)}{\hat{\omega}_2} \\
    &=& \mathcal{K}_{\gamma}(0, \tilde{\eta}_2(0)) \left[ \frac{1}{\tilde{\eta}_2(0)} \left.\frac{\partial \eta_2}{\partial \omega_2}\right|_{(0,\tilde{\eta}_2(0))} - \frac{1}{\hat{\omega}_2} \right].
 \end{eqnarray*}
 See (3.1.5) and (3.1.14) of \cite{Li:08} for more details.

\end{appendix}


\end{document}